\documentclass[a4wide]{article}
\usepackage[utf8]{inputenc}
\usepackage[a4paper, total={6.5in, 9.25in}]{geometry}

\usepackage{hyperref}
\usepackage[auth-lg]{authblk}

\usepackage[english]{babel}

\usepackage{amsmath}
\usepackage{graphicx}
\usepackage{mathtools}
\usepackage{amsfonts}
\usepackage{subfigure}
\usepackage{multirow}
\usepackage{array}
\newcolumntype{x}[1]{>{\centering\arraybackslash\hspace{0pt}}p{#1}}

\usepackage{color}

\usepackage[page]{appendix} 

\usepackage{amsthm}
\newtheorem*{remark}{Remark}

\begin{document}
\title{Uncertainty quantification for nonlinear solid mechanics using reduced order models with Gaussian process regression}

\author[1,2]{Ludovica Cicci} 
\author[1]{Stefania Fresca}
\author[3]{Mengwu Guo}
\author[1]{Andrea Manzoni}
\author[1]{Paolo Zunino}
\setlength{\affilsep}{1em}
\renewcommand\Authsep{, }
\affil[1]{MOX-Dipartimento di Matematica, Politecnico di Milano, Italy\thanks{\texttt{\{ludovica.cicci,stefania.fresca,andrea1.manzoni,paolo.zunino\}@polimi.it}}}
\affil[2]{School of Biomedical Engineering \& Imaging Sciences, King's College London, UK}
\affil[3]{Department of Applied Mathematics, University of Twente, the Netherlands\thanks{\texttt{m.guo@utwente.nl}}}
\date{}

\maketitle

\begin{abstract}
    Uncertainty quantification (UQ) tasks, such as sensitivity analysis and parameter estimation, entail a huge computational complexity when dealing with input-output maps involving the solution of nonlinear differential problems, because of the need to query expensive numerical solvers repeatedly. 
    Projection-based reduced order models (ROMs), such as the Galerkin-reduced basis (RB) method, have been extensively developed in the last decades to overcome the computational complexity of high fidelity full order models (FOMs), providing remarkable speedups when addressing UQ tasks related with parameterized  differential problems. Nonetheless, constructing a projection-based ROM that can be efficiently queried usually requires extensive modifications to the original code, a task which is non-trivial for nonlinear problems, or even not possible at all when proprietary software is used. Non-intrusive ROMs -- which rely on the FOM as a black box -- have been recently developed to overcome this issue. In this work, we consider ROMs exploiting proper orthogonal decomposition to construct a reduced basis from a set of FOM snapshots, and Gaussian process regression (GPR) to approximate the RB projection coefficients. Two different approaches, namely a global GPR and a tensor-decomposition-based GPR, are explored on a set of 3D time-dependent solid mechanics examples. Finally, the non-intrusive ROM is exploited to perform global sensitivity analysis (relying on both screening and  variance-based methods) and parameter estimation (through Markov chain Monte Carlo methods), showing remarkable computational speedups and very good accuracy compared to high-fidelity FOMs.
    
\end{abstract}

\renewcommand{\'}{``}


\section{Introduction}

Applied sciences and engineering problems, such as those arising in structural mechanics, are often described in terms of (time-dependent, nonlinear) partial differential equations (PDEs) that may be parameterized, i.e., involving several parameters to account for, e.g., different material properties, source terms, data, or geometrical features. 
Discrete, high-fidelity approximations of the PDE solutions can be computed by means of full order models (FOMs), which are usually computationally demanding - both in terms of CPU time and memory requirements - as fine computational grids and small time steps are needed to ensure accuracy of the solution. Despite the computational power available nowadays, relying on FOMs remains prohibitive in multi-query contexts, such as in the case of uncertainty quantification (UQ) tasks. 

When dealing with differential models hampered by uncertainty on the (possibly, many) input parameters, it may be useful for multiple purposes (i) to	detect and rank those inputs which need to be measured in order to reduce the output variance, (ii) to detect the parameters that have a better chance of being estimated in a subsequent estimation process, and (iii) to identify non-influential inputs in order to fix them to nominal values within their range of variability. To perform these tasks, global sensitivity analysis (SA) methods, such as screening methods \cite{Morris1991} and variance-based methods \cite{sobol1990sensitivity}, are commonly used in a variety of applications. However, computing sensitivity indices in the case of high-dimensional parameter spaces requires a large number of input-output evaluations (especially when SA is performed through Monte Carlo integration), so that the overall computational cost can become prohibitively expensive when relying on high-fidelity FOMs. Efficient surrogate models may be built \cite{le2017metamodel}, e.g., to compute the sensitivity indices directly using analytic formulas based on suitable regression predictors \cite{sudret2008global, cheng2020surrogate}, or to alleviate the cost of each forward simulation, eventually in a multi-fidelity framework \cite{qian2018multifidelity}.	

Another important problem addressed in real-world applications is the necessity to recover information about input parameters from limited and noisy observations of quantities of interest (QoIs), that are outputs computed using the problem solution. When addressing parameter estimation in a Bayesian framework, the solution to this inverse UQ problem is provided in terms of the posterior distribution of the unknown inputs conditioned on the observations. However, this distribution is not known analytically, and one needs to rely on sampling techniques, such as Markov chain Monte Carlo (MCMC), to characterize it and derive suitable statistics. Nonetheless, these methods require repeated (often, millions of) evaluations of the input-output map, so that relying on expensive high-fidelity FOMs becomes unfeasible. Different approaches can be employed to reduce the cost of solving parameter estimation problems \cite{frangos2010surrogate}, such as lowering the number of forward simulations required, reducing the dimensionality of the input space, or relying on efficient surrogate models for the solution to the forward problem \cite{li2014adaptive, Lassila2014, cui2015data, yang2021b}.

To tackle both SA and inverse UQ scenarios, we adopt a semi-intrusive approach by relying on a reduced order model (ROM) to compute efficient sensitivity estimates and posterior densities of the problem parameters, while keeping a good level of accuracy. ROMs featuring smaller dimension, i.e., a lower number of degrees of freedom (DOFs), have been extensively developed during the last decades to overcome the computational complexity of high-fidelity FOMs. Many of these reduction approaches rely on the assumption that the high-fidelity solution manifold, that is the set of all FOM solutions as the input parameters change, can be well approximated by a low-dimensional linear space, and employ the proper orthogonal decomposition (POD) to build such a reduced basis (RB) able to capture the most dominant features of the original system. The high-dimensional solution is thus reconstructed as a linear combination of these (hopefully) few modes with a small loss of accuracy. In particular, the evaluation  of high-fidelity solutions at several points in the time and parameter domains, the construction of the POD basis and the assembling of the ROM arrays, can be performed only once, offline. Then, for unseen time-parameter instances, the RB coefficients can be rapidly computed online, in order to reconstruct the high-dimensional solution. The key of success of these ROMs  lies in the full decoupling of the two stages, meaning that the online phase is independent of the high-fidelity dimension. ROMs can then be divided into intrusive and non-intrusive, according to the way the ROM approximation is generated and the reduced arrays are built.

Intrusive projection-based ROMs (see, e.g., \cite{benner2015survey} and references therein), such as reduced basis methods \cite{quarteroni2016reduced, hesthaven2016certified}, are constructed by projecting the high-fidelity model operators onto the linear subspace spanned by a set of problem-dependent basis functions. These methods have been successfully applied to a wide range of problems \cite{drohmann2012reduced, amsallem2012nonlinear, radermacher2016pod, ghavamian2017pod, bonomi2017reduced}, as they allow to retain the underlying structure of the high-fidelity FOM and to approximate the whole field variables that solve the FOM. Nonetheless, they require extensive code modifications, which might be non-trivial or even not possible when proprietary software is used. Moreover, in the case of nonlinear problems, or problems featuring a non-affine parametric dependence, assembling ROM operators requires a further level of hyper-reduction, relying on, e.g., the empirical interpolation method (EIM) \cite{barrault2004empirical} and its variants \cite{chaturantabut2010nonlinear, tiso2013modified, Lassila2014, negri2015efficient}. To overcome this difficulty, deep learning-based strategies have been recently developed \cite{Cicci2022}, showing a remarkable speed-up of the online phase also for complex time-dependent nonlinear applications \cite{cicci2022efficient} in the approximations of nonlinear ROM operators. Nonetheless, at every time step one must still solve the reduced problem and compute the solution vector, thus hampering the efficiency of the ROM in those cases where only scalar or low-dimensional QoIs need to be evaluated, or only the solution at prescribed time instances is required.

To overcome these issues, alternative data-driven methods can be used for the approximation of the RB coefficients without resorting to (Galerkin)  projection. In these cases, the FOM solution is projected onto the RB space and the combination coefficients are approximated by means of a surrogate model, exploiting, e.g., Gaussian process regression (GPR) \cite{guo2018reduced, guo2019, kast2020non, zhang2019model}, radial basis function (RBF) approximation \cite{audouze2013nonintrusive} or artificial neural networks (ANNs) \cite{hesthaven2018non, gao2021non, salvador20211}. The high-fidelity solver is thus used only offline as a 'black-box' to generate the necessary data for the RB construction and the training of the surrogate model. Non-intrusive POD-based approaches using RBFs to interpolate the ROM coefficients have been proposed in \cite{audouze2013nonintrusive, walton2013reduced,  xiao2017parameterized} and applied to nonlinear parameterized time-dependent PDEs. In \cite{hesthaven2018non} ANNs have been employed to build a regression model to compute the coefficients of a POD-based ROM, focusing on steady PDEs; a further extension to time-dependent nonlinear problems, i.e., unsteady flows, has been addressed in \cite{wang2019non}; see also, e.g., \cite{swischuk2019projection,yu2019non} for the use of machine learning strategies to approximate the POD coefficient. Alternatively, POD-GPR reduced order models have been developed in \cite{guo2018reduced} for steady nonlinear structural analysis and in \cite{guo2019} for time-dependent problems, and successfully used in engineering applications, such as electromagnetic scattering problems \cite{zhao2021non} and naval engineering problems \cite{ortali2020gaussian}. Non-intrusive ROMs using RBF, GPR, and ANN regressions are compared in \cite{berzicnvs2020standardized} and applied to complex flow problems in a time-dependent, deforming computational domain. See instead \cite{guo2022multi,conti2022multi} for an alternative use of ANNs to perform regression in the context of multi-fidelity methods relying on models of different fidelities. These latter may involve data-driven projection-based \cite{peherstorfer2016data, guo2022bayesian} ROMs or recently developed deep learning-based ROMs \cite{fresca2021comprehensive,fresca2021pod, franco2023deep, botteghi2022deep}.

In this work we apply non-intrusive ROMs for the efficient solution to nonlinear time-dependent parameterized problems arising in solid mechanics, in order to accelerate SA and the solution to Bayesian inverse problems. Once POD has been used to construct a low-dimensional linear subspace -- thus providing a way to embed physics into the ROM and to make it explainable --  the RB expansion coefficients are approximated by means of GPR models, taking into account (and comparing) two different approaches \cite{guo2018reduced} and \cite{guo2019}, both ensuring non-intrusiveness and online efficiency. With respect to other data-driven approaches, such as DNNs, the advantages of POD-GPR ROMs are twofold: on one hand, they require a relatively small amount of data to be trained; on the other hand, they allow to reconstruct the whole solution field, so that any quantity of interest can be evaluated without the need to retrain the model, whereas surrogate models would directly approximate the input-output map. Moreover, relying on GPRs allows to take into account the emulator uncertainty, which represents a further source of uncertainty in our problem.

Compared to the existing literature -- see, e.g., \cite{guo2018reduced, guo2019} -- we address in this paper the approximation of time-dependent, nonlinear problems set over three-dimensional domains, also featuring  a rather large number of input parameters. Moreover, a full UQ pipeline is shown, including both sensitivity analysis and parameter estimation, exploiting the non-intrusive ROM approximation enabled by the combined POD-GPR framework, and assessing its accuracy against the results obtained with a high-fidelity FOM based on the finite element method.

This work is structured as follows. In Section~\ref{sec:MOR} we review some basic facts about the RB method for time-dependent parameterized problems, describing in Section~\ref{sec:GPR} two different approaches based on GPRs for the efficient approximation of the RB coefficients. Numerical results related with two test cases in nonlinear solid mechanics are then presented in Section~\ref{sec:numerical_results}, where POD-GPR ROMs are compared in terms of efficiency and accuracy. In Section~\ref{sec:UQ} the solution to SA and inverse UQ problems is carried out by means of the POD-GPR ROMs. Finally, some conclusions are drawn in Section~\ref{sec:conclusions}.	
	

\section{RB methods for parametrized PDEs: basic facts and notation} \label{sec:MOR}

Our goal is to efficiently solve nonlinear time-dependent PDE problems depending on a set of input parameters $\boldsymbol{\mu}\in\mathcal{P}$, where $\mathcal{P} \subset \mathbb{R}^p$ is a compact set representing the parameter space, and $p$ is the number of inputs. In particular, we focus on problems governed by the elastodynamics equation
\begin{equation}\label{eq:state}
    \rho\partial_t^2 \boldsymbol{u}(\boldsymbol{X},t;\boldsymbol{\mu}) - \nabla_{\boldsymbol{X}}\cdot\boldsymbol{P}(\boldsymbol{u}(\boldsymbol{X},t;\boldsymbol{\mu});\boldsymbol{\mu}) = \mathbf{0} ~\text{in } \Omega_0\times\mathcal{T},
\end{equation}
suitably complemented by initial and boundary conditions. We denote the unknown state solution, i.e., the displacement field, by  $\boldsymbol{u}=\boldsymbol{u}(\boldsymbol{X},t;\boldsymbol{\mu})$,  depending on the spatial coordinate $\boldsymbol{X}$ in the reference configuration, the time variable by $t$, the parameter vector by $\boldsymbol{\mu}$, and the first Piola-Kirchhoff stress tensor by $\boldsymbol{P}$. Once this problem has been discretized in space relying, e.g., on the Galerkin finite element method, a dynamical system of dimension $N_h$ is obtained under the form
\begin{equation}\label{eq:disctete_state}
    \left\{\begin{array}{l}
        \rho_0\mathbf{M}\ddot{\mathbf{u}}_h(t;\boldsymbol{\mu}) + \mathbf{N}(\mathbf{u}_h(t;\boldsymbol{\mu}),t;\boldsymbol{\mu}) = \mathbf{F}^{ext}(t;\boldsymbol{\mu})\qquad \text{in}~\mathcal{T},\\
        \mathbf{u}_h(0;\boldsymbol{\mu})=\mathbf{u}_{h,0}(\boldsymbol{\mu}),\quad\dot{\mathbf{u}}_h(0;\boldsymbol{\mu})=\dot{\mathbf{u}}_{h,0}(\boldsymbol{\mu}),
    \end{array}\right.
\end{equation}
where $\mathbf{M}\in\mathbb{R}^{N_h\times N_h}$ denotes the mass matrix, $\mathbf{N}\in\mathbb{R}^{N_h}$ is the (nonlinear) vector of internal forces, $\mathbf{F}^{ext}\in\mathbb{R}^{N_h}$ is the vector of external forces and $\mathbf{u}_h(t)\in\mathbb{R}^{N_h}$ collects the degrees of freedom (DOFs) that represent the (semi-discrete) high fidelity solution at time $t$. Note that fully implicit schemes are employed for the sake of numerical stability, since they do not pose restrictions on the time step, otherwise required due to the highly nonlinear terms of the strain energy density function considered in this work. Here, $N_h$ denotes the number of DOFs, which depends on the underlying computational mesh and the space discretization used. Relying on finite differences for time discretization, and therefore approximating
\begin{align*}
    \dot{\mathbf{u}}_h(t^{n},\boldsymbol{\mu}) \approx \frac{\mathbf{u}_h(t^n,\boldsymbol{\mu})-\mathbf{u}_h(t^{n-1},\boldsymbol{\mu})}{\Delta t}, &&
    \ddot{\mathbf{u}}_h(t^n,\boldsymbol{\mu}) \approx \frac{\mathbf{u}_h(t^n,\boldsymbol{\mu})-2\mathbf{u}_h(t^{n-1},\boldsymbol{\mu})+\mathbf{u}_h(t^{n-2},\boldsymbol{\mu})}{\Delta t^2},
\end{align*}
where $\{t^0<t^1<\dots<t^{N_t}\}$ is a uniform partition of $\mathcal{T}$ with time step $\Delta t$, we obtain a sequence of fully-discrete problems that  in algebraic form read as follows: given $\boldsymbol{\mu}\in\mathcal{P}$, for $t^n$, $n=1,\dots, N_t$, find the discrete high-fidelity displacement $\mathbf{u}_h(t^n; \boldsymbol{\mu}) \in \mathbb{R}^{N_h}$ such that
\begin{equation} \label{eq:fully_discrete}
 \dfrac{\rho}{\Delta t^2}\mathbf{M}\mathbf{u}_h(t^{n},\boldsymbol\mu) + \mathbf{N}(\mathbf{u}_h(t^{n},\boldsymbol\mu),t^n;\boldsymbol{\mu}) - \dfrac{2\rho}{\Delta t^2}\mathbf{M}\mathbf{u}_h(t^{n-1},\boldsymbol\mu) + \dfrac{\rho}{\Delta t^2}\mathbf{M}\mathbf{u}_h(t^{n-2},\boldsymbol\mu) - \mathbf{F}^{ext}(t^{n},\boldsymbol\mu) = \mathbf{0},
\end{equation}
where $\mathbf{u}_h(t^{-1},\boldsymbol{\mu})$ and $\mathbf{u}_h(t^0,\boldsymbol{\mu})$ are given by the initial conditions. The solution to the nonlinear system (\ref{eq:fully_discrete}) can be finally computed by means of iterative algorithms, such as the Newton-Raphson method. However, to achieve sufficient accuracy for real-world applications, either the FOM dimension $N_h$ or the number of time steps $N_t$, as well as the total Newton iterations, can become extremely large, implying huge computational costs. 

Therefore, we rely on the RB method \cite{quarteroni2016reduced, hesthaven2016certified, hesthaven_pagliantini_rozza_2022} to seek a low-dimensional approximation of $\mathbf{u}_h(t;\boldsymbol{\mu})\in\mathbb{R}^{N_h}$ as a linear combination of $N\ll N_h$ parameter-independent basis functions constructed by performing proper orthogonal decomposition (POD) on a set of sampled FOM solutions. By denoting as $\mathbf{V}\in\mathbb{R}^{N_h\times N}$ the matrix collecting column-wise the DOFs of the RB functions, a reduced order solution $\mathbf{u}_h^{RB}(t;\boldsymbol{\mu})\in\mathbb{R}^{N_h}$ can be computed by projecting the FOM solution onto $\text{Col}(\mathbf{V})$, i.e., 
\begin{equation*}
    \mathbf{u}_h^{RB}(t;\boldsymbol{\mu})=\mathbf{VV}^T\mathbf{u}_h(t;\boldsymbol{\mu}),
\end{equation*}
so that we only need to recover the underlying map between $(t,\boldsymbol{\mu})$  and the reduced coefficients, that is
\begin{eqnarray*}
\mathbf{q}\colon(t,\boldsymbol{\mu})\mapsto\mathbf{V}^T\mathbf{u}_h(t;\boldsymbol{\mu}),
\end{eqnarray*}
for $(t;\boldsymbol{\mu})\in\mathcal{T}\times\mathcal{P}$. Note that, in this case, it holds $	\mathbf{u}_h^{RB}(t;\boldsymbol{\mu})= {\text{arg min}}_{\mathbf{x}^*_h\in\text{Col}(\mathbf{V})}\lVert \mathbf{u}_h(t;\boldsymbol{\mu})-\mathbf{x}^*_h\rVert$. However, when employing a Galerkin projection to compute the RB coefficients, the assembly of the reduced system still depends on the high-fidelity dimension due to nonlinearity of the problem, so that a further level of hyper-reduction is required to ensure computational efficiency, see, e.g., \cite{barrault2004empirical, chaturantabut2010nonlinear, negri2015efficient, tiso2013discrete}. The RB method thus yields an intrusive strategy to generate the reduced-order solution, with hyper-reduction techniques possibly hampering the convergence properties of the nonlinear solvers at the ROM level. Based on \cite{guo2018reduced, guo2019}, we resort to nonlinear regression to build efficient non-intrusive ROMs and obtain probabilistic distributions of the projection coefficients for each new value of the input vector, useful to perform UQ studies.
	
	
\section{Reduced order surrogate modeling with GP Regression} \label{sec:GPR}

To obtain the projection coefficients $\mathbf{q}(t;\boldsymbol{\mu})$, for $(t,\boldsymbol{\mu})\in\mathcal{T}\times\mathcal{P}$, and thus recover the RB solution $\mathbf{u}_h^{RB}(t;\boldsymbol{\mu})\approx\mathbf{u}_h(t;\boldsymbol{\mu})$, we rely on two different regression-based approaches, to which we refer to as \textit{global}, outlined in Section~\ref{sec:mono}, and \textit{tensor-decomposition-based}, proposed in \cite{guo2019} and described in Section~\ref{sec:segr}. We point out that, since the RB coefficients are uncorrelated, we rely on independent single-output GPs.

\subsection{Gaussian Process (GP) Regression}\label{sec:GP}

Given a dataset $\mathcal{D} = (\mathbf{X},\mathbf{y}) = \left\{ (\mathbf{x}^{(i)},y^{(i)})\in\mathcal{X}\times\mathcal{Y}\subset\mathbb{R}^{d}\times\mathbb{R}~|~ i=1,\dots,N_{data} \right\}$, where
\begin{equation*}
    \mathbf{X}=\bigl[\mathbf{x}^{(1)}|\dots|\mathbf{x}^{(N_{data})}\bigr] \quad\text{and}\quad \mathbf{y}=[y^{(1)},\dots,y^{(N_{data})}]^T
\end{equation*}
are the observed input matrix and the associated output vector, respectively, one wants to build a model $f\colon\mathcal{X}\rightarrow\mathcal{Y}$ which is able to infer the relationship between the $d$-dimensional independent variable $\mathbf{x}$ and the scalar dependent variable $y$, hence to determine the conditional distribution of the target output given the input. Since there may be more than a regression function $f$ that fit the data equally well, a Gaussian process (GP) \cite{williams2006gaussian} is used to assign a probability to each of these functions. Hence, given a mean $m(\mathbf{x})$ and a covariance function (or kernel) $\kappa(\mathbf{x},\mathbf{x}')$, one assumes that the prior of $f$ is a GP and write
\begin{equation*}
    f(\mathbf{x}) \sim \textup{GP}(m(\mathbf{x}),\kappa(\mathbf{x},\mathbf{x}')).
\end{equation*}

\begin{remark}
    Typically, the mean is chosen to be constant, e.g., equal to the mean of the training set (or to zero), whereas many different choices are available for the covariance $\kappa\colon\mathcal{X}\times\mathcal{X}\rightarrow\mathbb{R}$, such as the radial basis function (RBF) kernel
    \begin{equation}\label{eq:RBF}
        \kappa(\mathbf{x},\mathbf{x}') = \sigma_f^2~\textup{exp}\left( - \frac{\lVert\mathbf{x}-\mathbf{x}'\rVert^2}{2\ell^2} \right) = \sigma_f^2~\textup{exp}\left( -\frac{1}{2}\sum_{j=1}^d \frac{(x_j-x_j')^2}{\ell^2} \right),
    \end{equation}
    whose hyperparameters are the signal variance $\sigma_f^2$ and the length scale $\ell$. Alternatively, the automatic relevance determination (ARD) RBF kernel assumes each dimension as independent from the others,  thus allowing for different length scales $\ell_j$, for $j=1,\dots,d$. A non-stationary  covariance function may also be defined via a deep neural network architecture \cite{lee2017deep, guo2021brief}.
\end{remark}

The goal of GP regression is to incorporate the knowledge provided by the training data $\mathcal{D}$ into the prior to obtain the posterior (predictive) distribution. 
Even though in this work we rely on synthetic data, GPs can be trained on experimental outputs, so that measurement errors should be taken into account. Hence, in order to be as extensive as possible, we assume that the training output labels are corrupted by an additive independent identically distributed Gaussian noise $\epsilon\sim\mathcal{N}(0,\sigma_y^2)$, i.e., $y^{(i)} = f(\mathbf{x}^{(i)}) + \epsilon$, for $i=1,\dots,N_{data}$. Given a new set of testing point $\mathbf{x}^*\notin\mathcal{D}$, according to the theorem of conditional Gaussian (see, e.g., \cite{murphy2012machine}), the conditional distribution of the predicted GP realization is given by
\begin{equation*}
    f(\mathbf{x}^*)~|~\mathbf{X},\mathbf{y},\mathbf{x}^* \sim \textup{GP}(\widehat{m}(\mathbf{x}^*),\widehat{\kappa}(\mathbf{x}^*,{\mathbf{x}^*}')),
\end{equation*}
where
\begin{align*}
    \widehat{m}(\mathbf{x}^*) &= m(\mathbf{x}^*) + \kappa(\mathbf{x}^*,\mathbf{X})(\kappa(\mathbf{X},\mathbf{X})+\sigma_y^2\mathbf{I})^{-1}(\mathbf{y}-m(\mathbf{X})),\\
    \widehat{\kappa}(\mathbf{x}^*,{\mathbf{x}^*}') &= \kappa(\mathbf{x}^*,{\mathbf{x}^*}') - \kappa(\mathbf{x}^*,\mathbf{X})(\kappa(\mathbf{X},\mathbf{X})+\sigma_y^2\mathbf{I})^{-1}\kappa(\mathbf{X},\mathbf{x}^*).
\end{align*}

The hyperparameters of the GPR model, such as $\boldsymbol{\theta} = (\sigma_y,\sigma_f,\ell)$ for the RBF kernel, highly influence the predictive performances and should be optimized, e.g., by maximizing the log-likelihood
\begin{equation}\label{eq:optimizationGP}
    \log \pi(\mathbf{y}|\mathbf{X},\boldsymbol{\theta}) = -\frac{1}{2}\left(\mathbf{y}^T(\kappa(\mathbf{X},\mathbf{X})+\sigma_y^2\mathbf{I})^{-1}\mathbf{y} + \log\det(\kappa(\mathbf{X},\mathbf{X})+\sigma_y^2\mathbf{I})+N_{data}\log 2\pi\right).
\end{equation}
Once the hyperparameters have been optimized on the training set $\mathcal{D}$, the GPR model can be used to perform inference on an unseen test set, using for the regression the mean function calculated by the posterior, denoted from now on -- with a slight abuse of notation -- as $\widehat{f}$.

\begin{remark}\label{rmk:scaling}
    A common practice to handle inputs of varying magnitudes is to perform feature scaling on the training and test datasets. Widely used transformations of the data $\mathbf{v}\in\mathbb{R}^m$ and $\mathbf{M}\in\mathbb{R}^{m\times n}$ are:
    \begin{enumerate}
        \item the min-max normalization $\boldsymbol{\xi}^\mathbf{v}_\text{min-max}$ and $\boldsymbol{\xi}^\mathbf{M}_\text{min-max}$, such that for $i=1,\ldots, m$, 
        \begin{equation} \label{eq:min_max}
                (\boldsymbol{\xi}^\mathbf{v}_\text{min-max})_i\colon \mathbf{v}_{i} \mapsto \frac{\mathbf{v}_{i} - \min\mathbf{v}}{\max\mathbf{v}-\min\mathbf{v}}, \ \ \ (\boldsymbol{\xi}^\mathbf{M}_\text{min-max})_i \colon \mathbf{M}_{i:} \mapsto \frac{\mathbf{M}_{i:} - \min_k \mathbf{M}_{ik}}{\max_k\mathbf{M}_{ik}-\min_k\mathbf{M}_{ik}};
        \end{equation}  
        \item
        the standardization $\boldsymbol{\xi}^\mathbf{M}_\text{stand}$ or  $\boldsymbol{\xi}^\mathbf{v}_\text{stand}$, such that for $i=1,\ldots, m$, 
        \begin{equation} \label{eq:stand}
                (\boldsymbol{\xi}^\mathbf{v}_\text{stand})_i\colon \mathbf{v}_{i} \mapsto \frac{\mathbf{v}_{i} - \text{mean}~\mathbf{v}}{\text{std}~\mathbf{v}},  \ \ \ 
                (\boldsymbol{\xi}^\mathbf{M}_\text{stand})_i \colon \mathbf{M}_{i:} \mapsto \frac{\mathbf{M}_{i:} - \text{mean}_k~ \mathbf{M}_{ik}}{\text{std}_k~\mathbf{M}_{ik}}. 
        \end{equation} 
    \end{enumerate}
\end{remark}
	
	
\subsection{Global approach}\label{sec:mono}

In a first attempt, our goal is to directly learn the map between the time/parameter values and the projection coefficients onto the RB space, that is
\begin{equation*}
    (t,\boldsymbol{\mu}) \mapsto \mathbf{q}(t;\boldsymbol{\mu})= \mathbf{V}^T\mathbf{u}_h(t;\boldsymbol{\mu}),
\end{equation*}
by means of a regression model. In particular, we want to infer each RB coefficient $q_\ell(t;\boldsymbol{\mu})$, for $\ell=1,\dots,N$, by means of a GP trained on a $N_tN_s$-dimensional dataset of input-output pairs $(\mathbf{X},\mathbf{y}_\ell)$, where
\begin{equation*}
    \mathbf{X} = \left[\begin{array}{ccc|c|ccc}
        t^1 & \dots & t^{N_t} & \dots & t^1 & \dots & t^{N_t}\\
        \boldsymbol{\mu}_1 & \dots & \boldsymbol{\mu}_1 & \dots & \boldsymbol{\mu}_{N_s} & \dots & \boldsymbol{\mu}_{N_s}
    \end{array}\right]\in\mathbb{R}^{(p+1)\times N_tN_s},
\end{equation*}
\begin{equation*}
    \mathbf{y}_\ell = [q_\ell(t^1;\boldsymbol{\mu}_1),\dots,q_\ell(t^{N_t};\boldsymbol{\mu}_1)~|~\dots~|~q_\ell(t^1;\boldsymbol{\mu}_{N_s}),\dots,q_\ell(t^{N_t};\boldsymbol{\mu}_{N_s})]^T+\boldsymbol{\varepsilon}\in\mathbb{R}^{N_tN_s},
\end{equation*}
such that the mean of the posterior distribution, for $\ell=1,\dots,N$, can be used to make predictions for every new input $\mathbf{x}^*=(t^*,\boldsymbol{\mu}^*)\in\mathcal{T}\times\mathcal{P}$. This approach has been proposed in \cite{guo2018reduced} for nonlinear structural analysis, however restricted to quasi-static problems, and we refer to it as \textit{global}, since time and parameter are treated as a unique input variable. According to the notation of Section~\ref{sec:GP} we have $N_{data} = N_tN_s$, $\mathbf{x}^{(i)} = (t^{n_i},\boldsymbol{\mu}_{m_i})$ and $y^{(i)}=q_\ell(t^{n_i},\boldsymbol{\mu}_{m_i})+\varepsilon$, for $\ell=1,\dots,N$, where $n_i\in\{1,\dots,N_t\}$ and $m_i\in\{1,\dots,N_s\}$.

In this work we apply the transformations reported in Remark~\ref{rmk:scaling}, i.e., (\ref{eq:min_max}) or (\ref{eq:stand}), to $(\mathbf{X},\mathbf{X}^*)$ and $(\mathbf{y}_\ell,\widehat{f}(\mathbf{X}^*))$, for $\ell=1,\dots,N$, with $\mathbf{M}=\mathbf{X}$ and $\mathbf{v}=\mathbf{y}_\ell$, respectively, in order to handle values of different magnitudes. Note that the minimum, the maximum, the mean and the standard deviation are computed over the training arrays.
	

\subsection{Tensor-decomposition-based regression approach}\label{sec:segr}

A different approach is the so-called tensor-decomposition-based regression developed in \cite{guo2019}, here briefly described. Let $\mathbb{P}\in\mathbb{R}^{N_t\times N_s \times (p+1)}$ be such that
\begin{equation*}
    \mathbb{P}_{n,m,:} = \left[\begin{array}{c}
        t^n \\
        \boldsymbol{\mu}_m
    \end{array}\right], ~\text{for}~n=1,\dots,N_t,~m=1,\dots,N_s,
\end{equation*}
and let $\mathbb{Q}_\ell\in\mathbb{R}^{N_t\times N_s}$ be the corresponding output matrix for the $\ell$-th RB coefficient, that is
\begin{equation*}
    \mathbb{Q}_\ell = \left[\begin{array}{ccc}
        q_\ell(t^1;\boldsymbol{\mu}_1) & \dots & q_\ell(t^1;\boldsymbol{\mu}_{N_s})  \\
        \vdots & & \vdots \\
        q_\ell(t^{N_t};\boldsymbol{\mu}_1) & \dots & q_\ell(t^{N_t};\boldsymbol{\mu}_{N_s})
    \end{array}\right],
\end{equation*}
for $\ell = 1,\dots,N$. The (truncated) singular value decomposition (SVD) of $\mathbb{Q}_\ell$ can be written as
\begin{equation*}
    \mathbb{Q}_\ell\approx\boldsymbol{\Psi}^\ell\boldsymbol{\Lambda}^\ell(\boldsymbol{\Phi}^\ell)^T = \sum_{k=1}^{N_{q_\ell}} \lambda_k^\ell\boldsymbol{\Psi}_k^\ell(\boldsymbol{\Phi}_k^\ell)^T,
\end{equation*}
where $\boldsymbol{\Psi}_k^\ell$ and $\boldsymbol{\Phi}_k^\ell$ are the $k$-th discrete time- and parameter-modes, respectively, and the truncation rank $N_{q_\ell}$ is obtained for a chosen value of the SVD tolerance $\varepsilon_{SVD}>0$, according to the percentage of variance one wants to account for. By splitting $\mathbb{P}$ into two input matrices we define the following arrays
\begin{equation*}
    \mathbf{X}_t = \left[\begin{array}{c}
        t^1 \\
        \vdots \\
        t^{N_t}
    \end{array}\right]\in\mathbb{R}^{N_t\times 1}, \quad \mathbf{X}_{\boldsymbol{\mu}} = \left[\begin{array}{c}
        (\boldsymbol{\mu}_1)^T \\
        \vdots \\
        (\boldsymbol{\mu}_{N_s})^T
    \end{array}\right]\in\mathbb{R}^{N_s\times p}
\end{equation*}
related to the time instances and the model parameters, respectively, and the corresponding target outputs $\mathbf{y}_{\ell,t}:=\boldsymbol{\Psi}^\ell\in\mathbb{R}^{N_t\times N_{q_\ell}}$ and $\mathbf{y}_{\ell,\boldsymbol{\mu}}:=\boldsymbol{\Phi}^\ell\in\mathbb{R}^{N_s\times N_{q_\ell}}$, possibly corrupted by noise. 

For $k=1,\dots,N_{q_\ell}$, our aim is to the learn the $k$-th continuous time-mode $\widehat{\psi}_k^\ell(t)$ given $\left\{\mathbf{X}_t,(\mathbf{y}_{\ell,t})_{1:N_t,k}\right\}$, and the $k$-th continuous parameter-mode $\widehat{\phi}_k^\ell(\boldsymbol{\mu})$ given $\left\{\mathbf{X}_{\boldsymbol{\mu}},(\mathbf{y}_{\ell,\boldsymbol{\mu}})_{1:N_s,k}\right\}$, so that, for a new value $\mathbf{x}^*=(t^*;\boldsymbol{\mu}^*)\in\mathcal{T}\times\mathcal{P}$, we can compute
\begin{equation*}
    q_\ell(\mathbf{x}^*) \approx \widehat{q}_\ell(\mathbf{x}^*) = \sum_{k=1}^{N_{q_\ell}} \lambda_k^\ell \widehat{\psi}_k^\ell(t^*) \widehat{\phi}_k^\ell(\boldsymbol{\mu}^*).
\end{equation*}

\begin{remark}\label{rmk:TD_number_GPR}
    For each component of the truncated SVD, i.e., for $k=1,\dots,N_{q_\ell}$, two Gaussian processes -- one for the time variable and the other for the parameter vector variable -- are trained, so that $2\sum_{\ell=1,\dots,N}N_{q_\ell}$ GPs are considered. Nonetheless, the dimensions of the training datasets are now much smaller with respect to the previous global approach ($N_t$ or $N_s$ instead of $N_tN_s$), thus reducing the overall computational time required for computing the optimal hyper-parameters.
\end{remark}

As before, we employ two scaling strategies, namely the min-max normalization and the standardization. In this case, the transformations in (\ref{eq:min_max}) or in (\ref{eq:stand}) are applied to the time- and parameter-datasets separately. To be more specific, the input datasets $(\mathbf{X}_t,\mathbf{X}^*_t)$ and $(\mathbf{X}_{\boldsymbol{\mu}},\mathbf{X}^*_{\boldsymbol{\mu}})$ are scaled according to maps $\boldsymbol{\xi}^{\mathbf{X}_t}_\star$ and $\boldsymbol{\xi}^{\mathbf{X}_{\boldsymbol{\mu}}}_\star$, respectively, where $\star\in\{\text{min-max},\text{stand}\}$. Meanwhile, for $\ell=1,\dots,N$ and $k=1,\dots,N_{q_\ell}$, the output datasets are scaled by applying $\boldsymbol{\xi}^{(\mathbf{y}_{\ell,t})_{1:N_t,k}}_\star$ to $((\mathbf{y}_{\ell,t})_{1:N_t,k},\widehat{\psi}_k^\ell(\mathbf{X}^*_t))$ and $\boldsymbol{\xi}^{(\mathbf{y}_{\ell,\boldsymbol{\mu}})_{1:N_s,k}}_\star$ to $((\mathbf{y}_{\ell,\boldsymbol{\mu}})_{1:N_s,k},\widehat{\phi}_k^\ell(\mathbf{X}^*_{\boldsymbol{\mu}}))$.


\section{Numerical results}\label{sec:numerical_results}

In this section we present the results obtained on two benchmark problems inspired by those used for verification of cardiac mechanics software \cite{land2015verification}, namely the deformation of a clamped rectangular beam and the active contraction of a truncated ellipsoid, governed by the elastodynamics equation (\ref{eq:state}). In particular, we focus on soft tissue mechanics \cite{holzapfel2003biomechanics}, a vivid area of research in computational mechanics characterized by many difficulties due do the presence of strong material anisotropy, large deformations, highly nonlinear stress-strain behaviors. Although specific for this field, the proposed test cases can also be of interest for other fields of application sharing similar features. 
For both test cases, we evaluate the performances of the global and the tensor-decomposition-based approaches, either in terms of training time, computational speed-up with respect to the high fidelity model and accuracy. Specifically, to test the fidelity of the GPR predictions on the RB coefficients $q_\ell$, for $\ell=1,\dots,N$, over the test set, we use the mean square error (MSE)
	\begin{equation*}
		\text{MSE}(q_\ell) = \frac{1}{N_s^{test}N_t^{test}}\sum_{m=1}^{N_s^{test}}\sum_{n=1}^{N_t^{test}} (q_\ell(t^n;\boldsymbol{\mu}_m)-\widehat{q}_\ell(t^n;\boldsymbol{\mu}_m))^2
	\end{equation*}
	and the relative squared error (RSE)
	\begin{equation*}
		\text{RSE}(q_\ell) = \frac{\sum_{m=1}^{N_s^{test}}\sum_{n=1}^{N_t^{test}} (q_\ell(t^n;\boldsymbol{\mu}_m)-\widehat{q}_\ell(t^n;\boldsymbol{\mu}_m))^2}{\sum_{m=1}^{N_s^{test}}\sum_{n=1}^{N_t^{test}} (q_\ell(t^n;\boldsymbol{\mu}_m)-\bar{q}_\ell)^2}, ~\text{with}~\bar{q}_\ell = \frac{1}{N_s^{test}N_t^{test}} \sum_{m=1}^{N_s^{test}}\sum_{n=1}^{N_t^{test}} q_\ell(t^n;\boldsymbol{\mu}_m).
	\end{equation*}
	The overall accuracy of the ROMs with respect to the FOM is evaluated using the time-averaged absolute and relative $L^2$-errors defined as
	\begin{align*}
		\text{tAE}(\boldsymbol{\mu}) &= \frac{1}{N_t^{test}}\sum_{n=1}^{N_t^{test}}\lVert \mathbf{u}_h(t^n;\boldsymbol{\mu})- \mathbf{V}\widehat{\mathbf{q}}(t^n;\boldsymbol{\mu})\rVert_{L^2}, \\
		\text{tRE}(\boldsymbol{\mu}) &= \frac{1}{N_t^{test}}\sum_{n=1}^{N_t^{test}}\frac{\lVert \mathbf{u}_h(t^n;\boldsymbol{\mu})- \mathbf{V}\widehat{\mathbf{q}}(t^n;\boldsymbol{\mu})\rVert_{L^2}}{\lVert\mathbf{u}_h(t^n;\boldsymbol{\mu})\rVert_{L^2}},
	\end{align*}
	respectively, and their mean over the test set, i.e., $\underset{m=1,\dots,N_s^{test}}{\text{mean}}~\text{tAE}(\boldsymbol{\mu}_m)$ and  $\underset{m=1,\dots,N_s^{test}}{\text{mean}}~\text{tRE}(\boldsymbol{\mu}_m)$. 
	
	The GPR models are implemented in Python using the library GPy \cite{gpy2014}, a Gaussian processes framework. All the computations have been performed on a Notebook with 2.60GHz Intel Core i7-9750H CPU and 16GB RAM. 
	
	\begin{remark}
		Numerical simulations have been run in serial and GPs have been trained and tested sequentially, i.e., on a outer loop from $\ell=1$ to $\ell=N$. However, better performances with respect to CPU time can be achieved by exploiting the independence of the GPR models associated to different RB coefficients, thus performing training and testing in parallel.	
	\end{remark}

\subsection{Benchmark tests for nonlinear solid mechanics}
 
For the setup of the benchmark test cases, we partition the boundary in three regions, corresponding to displacement, pressure and stress-free conditions respectively. More precisely, we set  Dirichlet boundary conditions on a specific portion of the boundary $\partial\Omega_0^{D}$ and the pressure load
\begin{equation*}
    \boldsymbol{p}(t;\boldsymbol{\mu}) = \widetilde{p}\frac{t}{T}
\end{equation*}
is applied on $\partial\Omega_0^{P}$. Moreover, Neumann conditions are set on all other surfaces, i.e., $\partial\Omega_0^{N} = \partial\Omega_0\backslash(\partial\Omega_0^{D}\cup\partial\Omega_0^{P})$, and homogeneous initial conditions are prescribed at time $t=0$~s. 

Equation (\ref{eq:state}) must be further complemented with a suitable constitutive relation for the material which can be expressed, e.g., as a nonlinear operator mapping the Green-Lagrange strain tensor $$\boldsymbol{E}(\boldsymbol{u}) = \frac{1}{2}(\boldsymbol{F}(\boldsymbol{u})^T\boldsymbol{F}(\boldsymbol{u})-\boldsymbol{I})$$ to the first Piola-Kirchhoff stress tensor $\boldsymbol{P}(\boldsymbol{u})$, where $\boldsymbol{F}(\boldsymbol{u})=\boldsymbol{I}+\nabla_{\boldsymbol{X}}\boldsymbol{u}$ is the deformation gradient. In this work we consider a nearly-incompressible hyperelastic material with the following strain-energy function
\begin{equation}\label{eq:guccione_energy}
    \mathcal{W}(\boldsymbol{F}) = \frac{C}{2}(e^{Q(\boldsymbol{E})-1}) + \frac{K}{2}(J-1)\ln(J),
\end{equation}
such that $\boldsymbol{P}(\boldsymbol{u}) = \partial\mathcal{W}(\boldsymbol{F})/\partial \boldsymbol{F}$. The first term is given by the constitutive law proposed in \cite{guccione1995finite}, and the second one is used to enforce incompressibility. Here, $C>0$ is a scaling constant, $J=\det(\boldsymbol{F}(\boldsymbol{u}))$ measures the change in volume during motion, and $K>0$ is the bulk modulus penalizing large volume variations. In particular, we consider the following form for $Q$,
\begin{equation*}
    Q = b_xE_{xx}^2 + b_yE_{yy}^2 + b_zE_{zz}^2 + b_{xy}(E_{xy}^2 + E_{yx}^2) + b_{xz}(E_{xz}^2 + E_{zx}^2) + b_{yz}(E_{yz}^2 + E_{zy}^2),
\end{equation*}
where $b_i$, for $i\in\{x,y,z\}$, are related to the material stiffness in the different directions, and $E_{ij}$, for $i,j\in\{x, y, z\}$, are the components of the strain tensor $\boldsymbol{E}(\boldsymbol{u})$. Therefore, we end up with the following initial-boundary value problem: given $\boldsymbol{\mu}\in\mathcal{P}$, find the parameterized displacement field $\boldsymbol{u}(\boldsymbol{\mu})\colon\Omega_0\times(0,T)\rightarrow\mathbb{R}^3$ such that
\begin{equation}\label{eq:PDE}
    \left\{\begin{array}{lll}
        \rho\partial_t^2 \boldsymbol{u}(\boldsymbol{X},t;\boldsymbol{\mu}) - \nabla_{\boldsymbol{X}}\cdot\boldsymbol{P}(\boldsymbol{u}(\boldsymbol{X},t;\boldsymbol{\mu});\boldsymbol{\mu}) = \mathbf{0} & \text{in} & \Omega_0\times(0,T)\\
        \boldsymbol{u}(\boldsymbol{X},t;\boldsymbol{\mu}) = \mathbf{0} & \text{on} & \partial\Omega_0^{D}\times(0,T)\\
        \boldsymbol{P}(\boldsymbol{u}(\boldsymbol{X},t;\boldsymbol{\mu});\boldsymbol{\mu})\boldsymbol{\nu} = -\boldsymbol{p}(t;\boldsymbol{\mu})J\boldsymbol{F}(\boldsymbol{u}(\boldsymbol{X},t;\boldsymbol{\mu});\boldsymbol{\mu})^{-T}\nu & \text{on} & \partial\Omega_0^{P}\times(0,T)\\
        \boldsymbol{P}(\boldsymbol{u}(\boldsymbol{X},t;\boldsymbol{\mu});\boldsymbol{\mu})\boldsymbol{\nu} = \mathbf{0} & \text{on} & \partial\Omega_0^{N}\times(0,T)\\
        \boldsymbol{u}(\boldsymbol{X},t;\boldsymbol{\mu}) = \mathbf{0}, \partial_t\boldsymbol{u}(\boldsymbol{X},t;\boldsymbol{\mu}) = \mathbf{0} & \text{in} & \Omega_0\times\{0\},
    \end{array}\right.
\end{equation}
where $\boldsymbol{\nu}$ denotes the outer normal unit vector.	
	
\subsubsection{Test 1: deformation of a beam} \label{sec:beam}

\begin{figure}[h!]
    \centering
    \includegraphics[width=\textwidth]{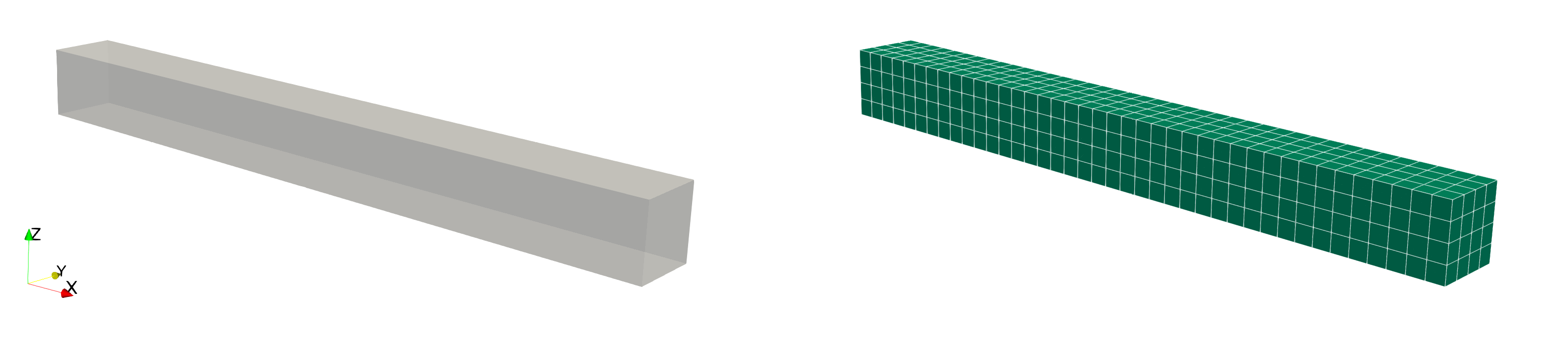}
    \caption{Deformation of a beam. Reference geometry (left) and computational mesh (right).}
    \label{fig:beam}
\end{figure}

The first benchmark test deals with the deformation of a clamped rectangular beam due to a pressure load acting on the bottom face ot the body. The problem is defined on the reference domain $\Omega_0 = [0,10^{-2}]\text{m}\times[0,10^{-3}]\text{m}\times[0,10^{-3}]\text{m}$, reported in Figure~\ref{fig:beam} together with the computational hexahedral mesh. With respect to system (\ref{eq:PDE}), we define $\partial\Omega_0^{D}=\Omega_0\cap\{x=0\}$ and $\partial\Omega_0^{P}=\Omega_0\cap\{z=0\}$.

As input parameters we consider those of the Guccione's law, i.e., $b_f,b_s,b_n,b_{fs},b_{fn},b_{sn},K$ and $G$, and the slope $\widetilde{p}$ of the pressure load, such that 
    $\boldsymbol{\mu}=[b_f,b_s,b_n,b_{fs},b_{fn},b_{sn},K,G,\widetilde{p}]\in\mathcal{P}
= \bigcup_{i=1}^9[l_i,r_i] \subset\mathbb{R}^9$, see Table~\ref{tab:beam_parameter_space}.	

\begin{table}[h!]
    \centering
    \begin{tabular}{|c|ccccccccc|}
        \hline
        & $b_f$ & $b_s$ & $b_n$ & $b_{fs}$ & $b_{fn}$ & $b_{sn}$ & $K$ & $C$ & $\widetilde{p}$ \\
        & & & & & & & [kPa] & [kPa] & [kPa]  \\
        \hline
        $l_i$ & $ 4$ & $1$ & $1$ & $2$ & $2$ & $1$ & $25$ & $1$ & $0.002$ \\
        $r_i$ & $12$ & $3$ & $3$ & $6$ & $6$ & $3$ & $75$ & $3$ & $0.006$ \\
        \hline
    \end{tabular}
    \caption{Deformation of a beam. Lower and upper bounds for the components of the parameter vector.}
    \label{tab:beam_parameter_space}
\end{table}

The FOM is built on an hexahedral mesh with 1025 vertices using $\mathbb{Q}_1$-FE in space and a BDF1 scheme in time with time step $\Delta t = 0.005$~s and final time $T=0.25$~s, so that $N_h=3075$ and $N_t=50$. In this setting the average CPU time required for computing the FOM solution for a given parameter vector is $48$~s.

With the aim of analyzing and comparing the two POD-GPR approaches, namely the global the and tensor-decomposition-based, we collect full-order snapshots at $N_s=100$ parameter locations sampled in $\mathcal{P}$ by latin hypercube sampling (LHS) for the construction of the RB functions. The same set solution vectors is then split (8:2) into training and testing datasets, and used for the construction of the GPR models. We point out that, unless otherwise specified, we use the same time points for training and testing, so that $N_t^{test} = N_t$.

\begin{remark}
    For the POD-based ROMs to accurately approximate the FOM displacement field for any input parameter, one needs to guarantee that the snapshots used for the construction of the RB basis sufficiently capture the solution dynamics, so that the choice of the sampled parameter locations highly influences the accuracy of the ROM.
    Relying on random strategies may require a large number of high-fidelity snapshots, leading to extremely expensive offline computational times. To overcome these issues, suitable adaptive sampling strategies have been adopted, e.g., in \cite{guo2018reduced, alsayyari2019nonintrusive, gao2021non, chen2018greedy}.
\end{remark}	
	
\subsubsection{Test 2: active contraction of a truncated ellipsoid}\label{sec:prolate}

\begin{figure}[h]
    \centering
    \includegraphics[width=0.75\textwidth]{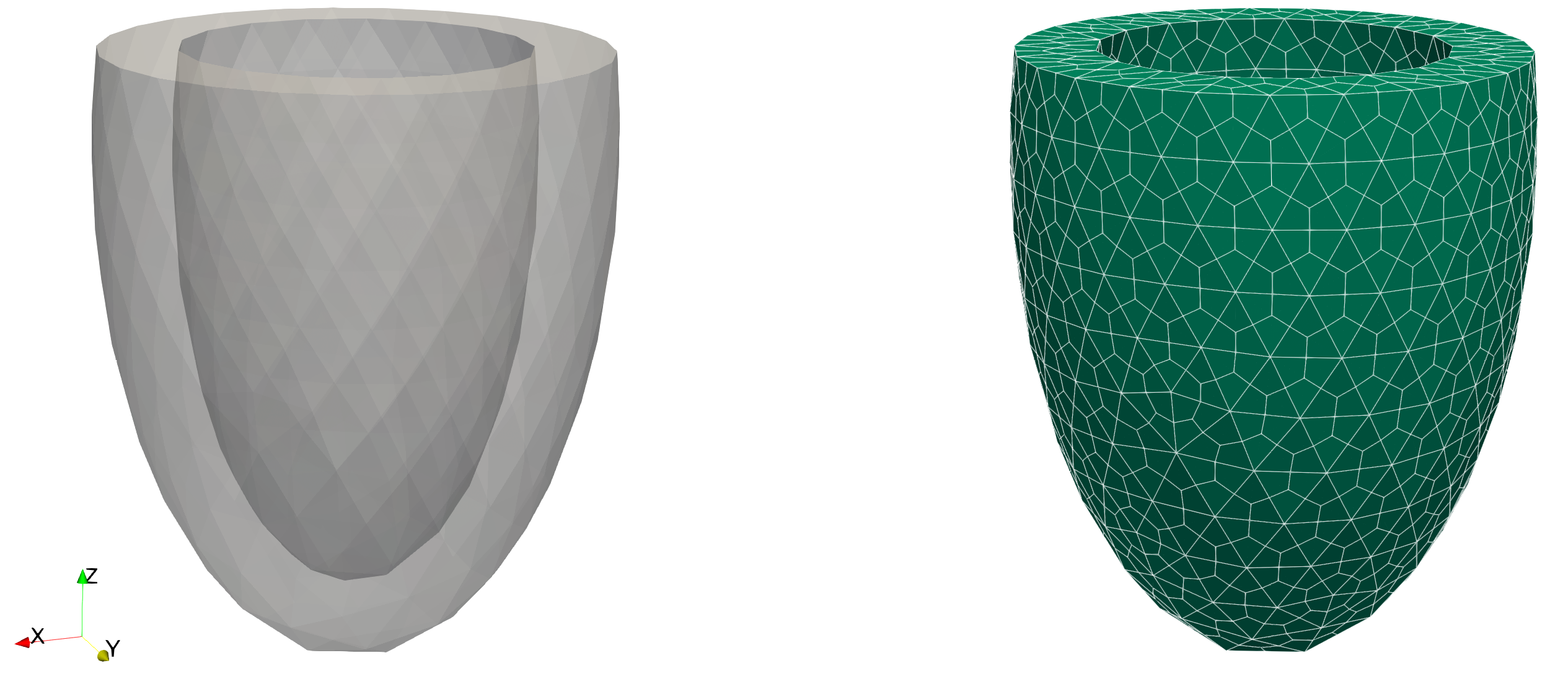}
    \caption{Active contraction of a truncated ellipsoid. Reference geometry (left) and computational mesh (right) of a truncated ellipsoid.}
    \label{fig:prolate}
\end{figure}

For the second problem we assume as reference geometry $\Omega_0$ the truncated ellipsoid described in \cite{land2015verification} and reported in Figure~\ref{fig:prolate}, and consider a nearly-incompressible material governed by the \textit{Guccione et al.} strain-energy function (\ref{eq:guccione_energy}). Moreover, we integrate active tension into the passive stress tensor by adding a parameterized time-dependent term which is assumed to act only in the fiber direction, that is
\begin{equation*}
    \mathbf{P}(\mathbf{u}) = \frac{\partial\mathcal{W}(\mathbf{F})}{\partial\mathbf{F}} + T_a(t;\boldsymbol{\mu})(\mathbf{Ff}\otimes\mathbf{f}),
\end{equation*}
where $\mathbf{f}\in\mathbb{R}^3$ denotes the reference unit vector in the fiber direction and we assume
\begin{equation*}
    T_a(t;\boldsymbol{\mu}) = \widetilde{T}_a\frac{t}{T},
\end{equation*}
for $\widetilde{T}_a>0$. The base plane $\partial\Omega_0^{D}=\Omega_0\cap\{z=\bar{z}\}$ is fixed in all directions, and an external pressure in applied at the inner boundary $\partial\Omega_0^{P}$. Neumann conditions are finally applied at the outer surface $\partial\Omega_0^{N}$.

In addition to the material parameters and the slope $\widetilde{p}$ of the pressure load, we consider the maximum value $\widetilde{T}_a$ of the active tension  and the fiber orientation angles $\boldsymbol{\alpha}_{epi}$ and $\boldsymbol{\alpha}_{endo}$ at at the epicardium and the endocardium, respectively. Hence, $\boldsymbol{\mu} =[b_f,b_s,b_n,b_{fs},b_{fn},b_{sn},K,G,\widetilde{p},\widetilde{T}_a,\boldsymbol{\alpha}_{epi},\boldsymbol{\alpha}_{endo}]\in\mathcal{P}\subset\mathbb{R}^{12}$, where the parameter space $\mathcal{P}$ is reported in Table~\ref{tab:prolate_parameter_space}.

\begin{table}[h!]
    \centering
    \begin{tabular}{|c|cccccccccccc|}
        \hline
        & $b_f$ & $b_s$ & $b_n$ & $b_{fs}$ & $b_{fn}$ & $b_{sn}$ & $K$ & $C$ & $\widetilde{p}$ & $\widetilde{T}_a$ & $\boldsymbol{\alpha}_{epi}$ & $\boldsymbol{\alpha}_{endo}$ \\
        & & & & & & & [kPa] & [kPa] & [kPa] & [kPa] & [$^\circ$]& [$^\circ$] \\
        \hline
        $l_i$ & $6.6$ & $1.65$ & $1.65$ & $3.3$ & $3.3$ & $1.65$ & $40$ & $1.5$ & $14$ & $49.5$ & $-105.5$ & $74.5$ \\
        $r_i$ & $9.4$ & $2.35$ & $2.35$ & $4.7$ & $4.7$ & $2.35$ & $60$ & $2.5$ & $16$ & $70.5$ & $-74.5$  & $105.5$\\
        \hline
    \end{tabular}
    \caption{Active contraction of a truncated ellipsoid. Parameters' range.}
    \label{tab:prolate_parameter_space}
\end{table}

The FOM is built on an hexahedral mesh with $6455$ vertices using $\mathbb{Q}_1$-FE, thus obtaining $N_h = 19365$, and a BDF1 scheme in time with $\Delta t = 0.005$~s and $T=0.25$~s, such that $N_t = 50$. For the construction and evaluation of the reduced model, we sample $100$ points in $\mathcal{P}$ using LHS and compute the corresponding high-fidelity displacements, which are then split (8:2) into training and test sets. In this case, the CPU time required to compute a single FOM solution is in average $6$ min $11$ s.

\subsection{Comparison of global and tensor-decomposition POD-GPR approaches}

\subsubsection{Test 1: deformation of a beam} \label{sec:beam_results}

Reduced basis of different dimensions $N$ can be constructed by performing POD on the FOM snapshots matrix $\mathbf{S}_d$ using different truncation tolerances $\varepsilon_{POD}$, as reported in Table~\ref{tab:beam_POD}. For the case at hand, 
$N=5$ RB functions (corresponding to a POD tolerance of $5\cdot10^{-4}$) are sufficient to obtain a good approximation of the high-fidelity FOM solutions, given the fast decay of the singular values.

\begin{table}[h!]
    \centering
    \begin{tabular}{|c|ccccc|}
        \hline
        $\varepsilon_{POD}$ & $10^{-3}$ & $5\cdot10^{-4}$ & $10^{-4}$ & $5\cdot10^{-5}$ & $10^{-5}$\\
        $N$ & 4 & 5 & 8 & 12 & 22 \\
        \hline
    \end{tabular}
    \caption{Deformation of a beam. Dimension of RB subspaces for different POD tolerances.}
    \label{tab:beam_POD}
\end{table}
	
Since the GPR models highly depend on the choice of the covariance functions, we analyze their predictive performances using different kernels,  namely polynomial, RBF and ARD-RBF, and different scaling techniques, i.e., standardization and min-max normalization. For this comparison we consider only the first $N_s=50$ training samples collected offline for the construction of the RB matrix $\mathbf{V}$. We point out that, for the tensor-decomposition-based approach, time- and parameter-modes are trained as GPs after the truncated SVD with tolerance $10^{-2}$ is computed, obtaining $3\leq N_\ell\leq6$. 
We present here the results subdivided in two parts: (i) first we address the GPRs predictions on the RB coefficients; (ii) then we assess the overall accuracy of the ROM with respect to the FOM.

\begin{figure}[h!]
    \centering
    \includegraphics[width=0.95\textwidth]{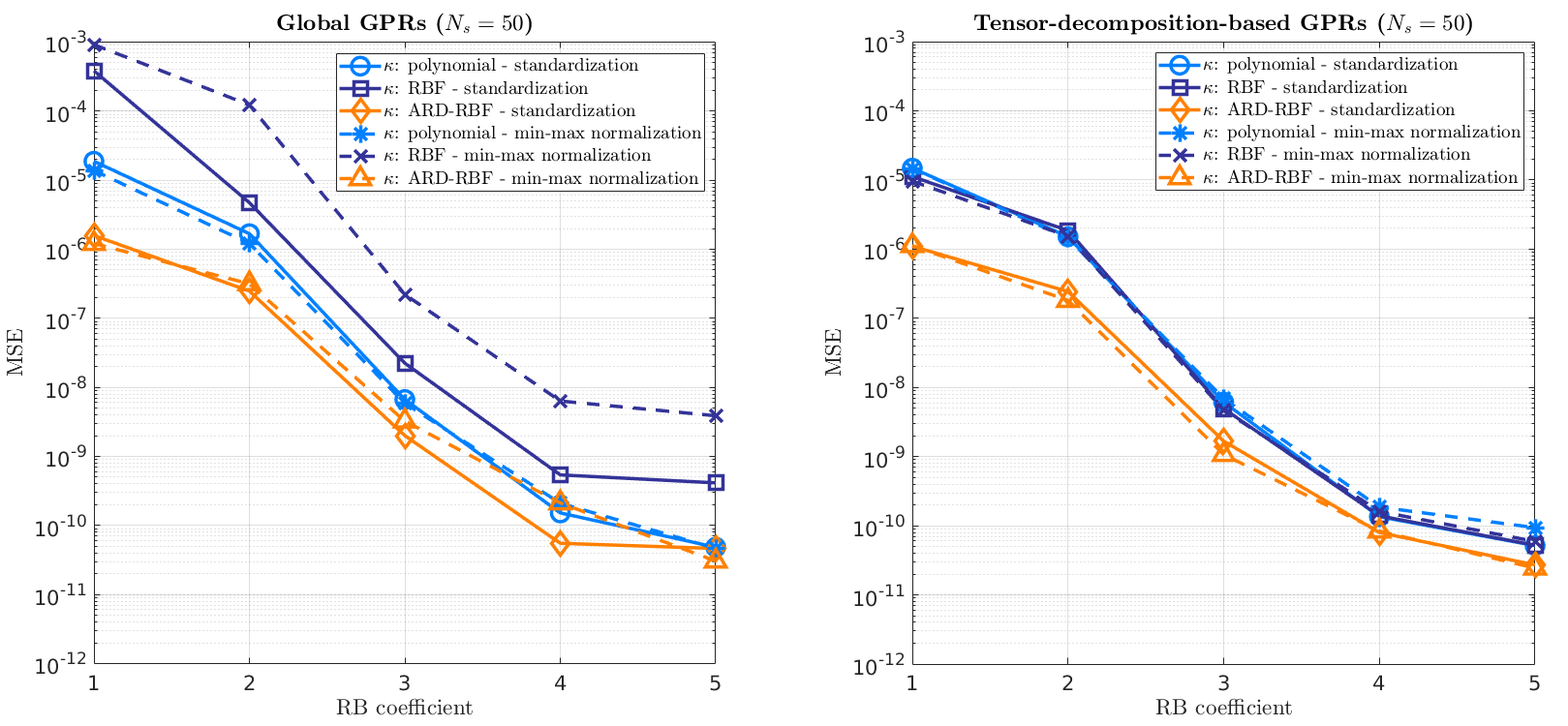}
    \caption{Deformation of a beam. Mean square error of global (left) and tensor-decomposition-based (right) POD-GPR ROMs computed using different kernel types and scaling techniques, for $N=5$.}
    \label{fig:beam_Ns_kernel_scaling}
\end{figure}

 \paragraph{Predictions of GPRs on RB coefficients.}
Figure~\ref{fig:beam_Ns_kernel_scaling} reports the MSEs computed using different covariance functions and scaling techniques, for both global and tensor-decomposition-based GPR models. In this case, the ARD-RBF kernel is more accurate than the other covariance functions, whereas the standardization and min-max normalization are almost comparable, with slightly better accuracy obtained using the former option. As a consequence, from now on we will focus only on ARD-RBF kernels and scale the dataset according to (\ref{eq:stand}).
	
\begin{figure}[h!]
    \centering
    \includegraphics[width=0.95\textwidth]{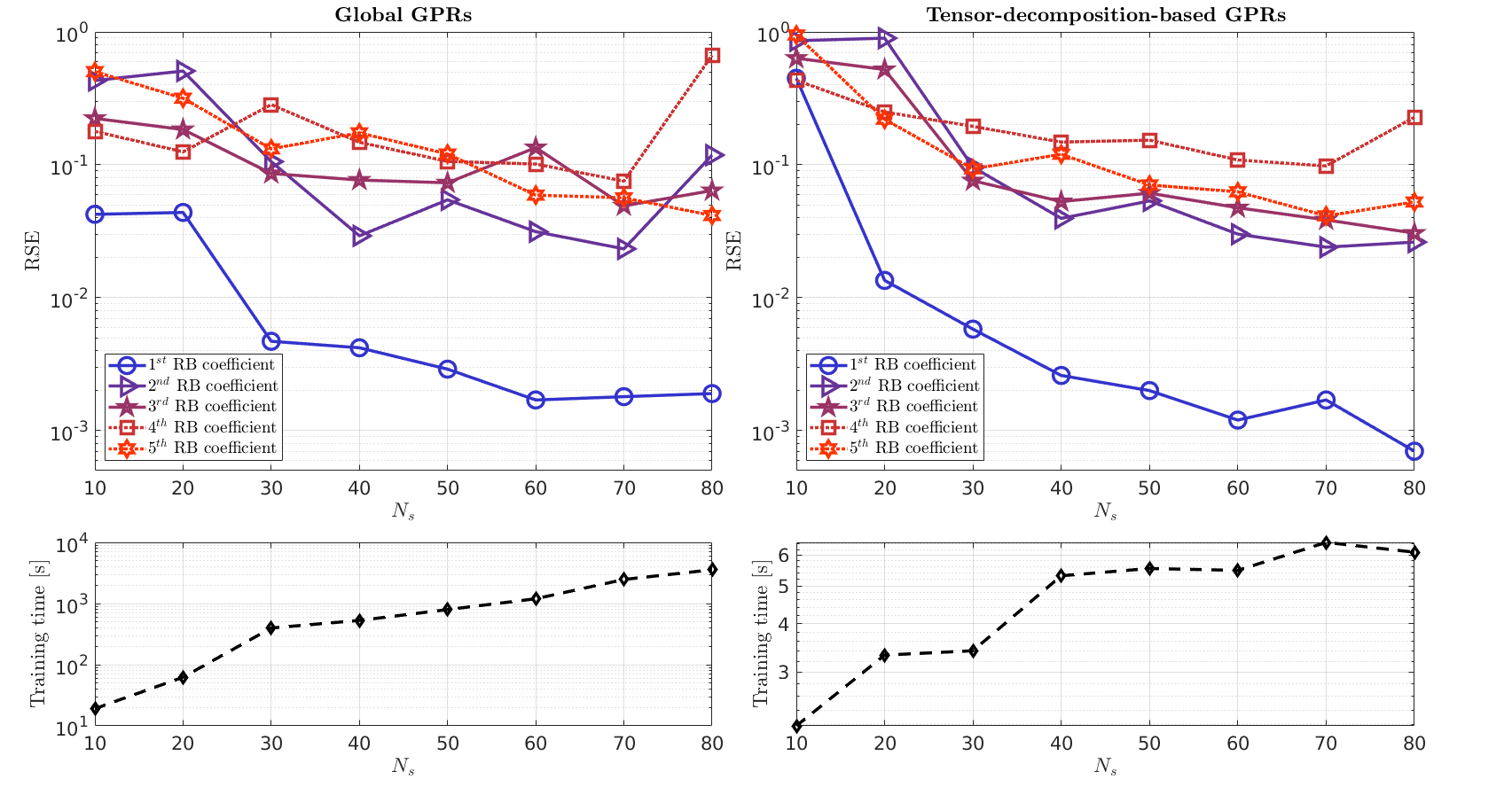}
    \caption{Deformation of a beam. Relative squared error of the global (top left) and the tensor-decomposition-based (top right) POD-GPR ROMs over the size of the training set. Corresponding training times (bottom) are also reported. These results have been obtained taking into account a testing set of $N_s^{test}=20$ parameters, unseen during training.}\label{fig:beam_convergence_mono_segr}
\end{figure}

We analyze the performances of the global and the tensor-decomposition-based POD-GPR ROMs in terms of accuracy and efficiency by increasing the size of the training set from $N_s=10$   to $N_s=80$. The RSE
computed over the testing set is reported in Figure~\ref{fig:beam_convergence_mono_segr}, showing a good improvement of the approximation  by both approaches as we increase the size of the training set from $N_s=10$ to $N_s=40$. We point out that the tensor-decomposition-based approach beats the global one in terms of computational efficiency during training, as this is performed for each GP on a much smaller dataset (see Remark~\ref{rmk:TD_number_GPR}). 

\begin{figure}[h!]
    \centering
    \includegraphics[width=0.95\textwidth]{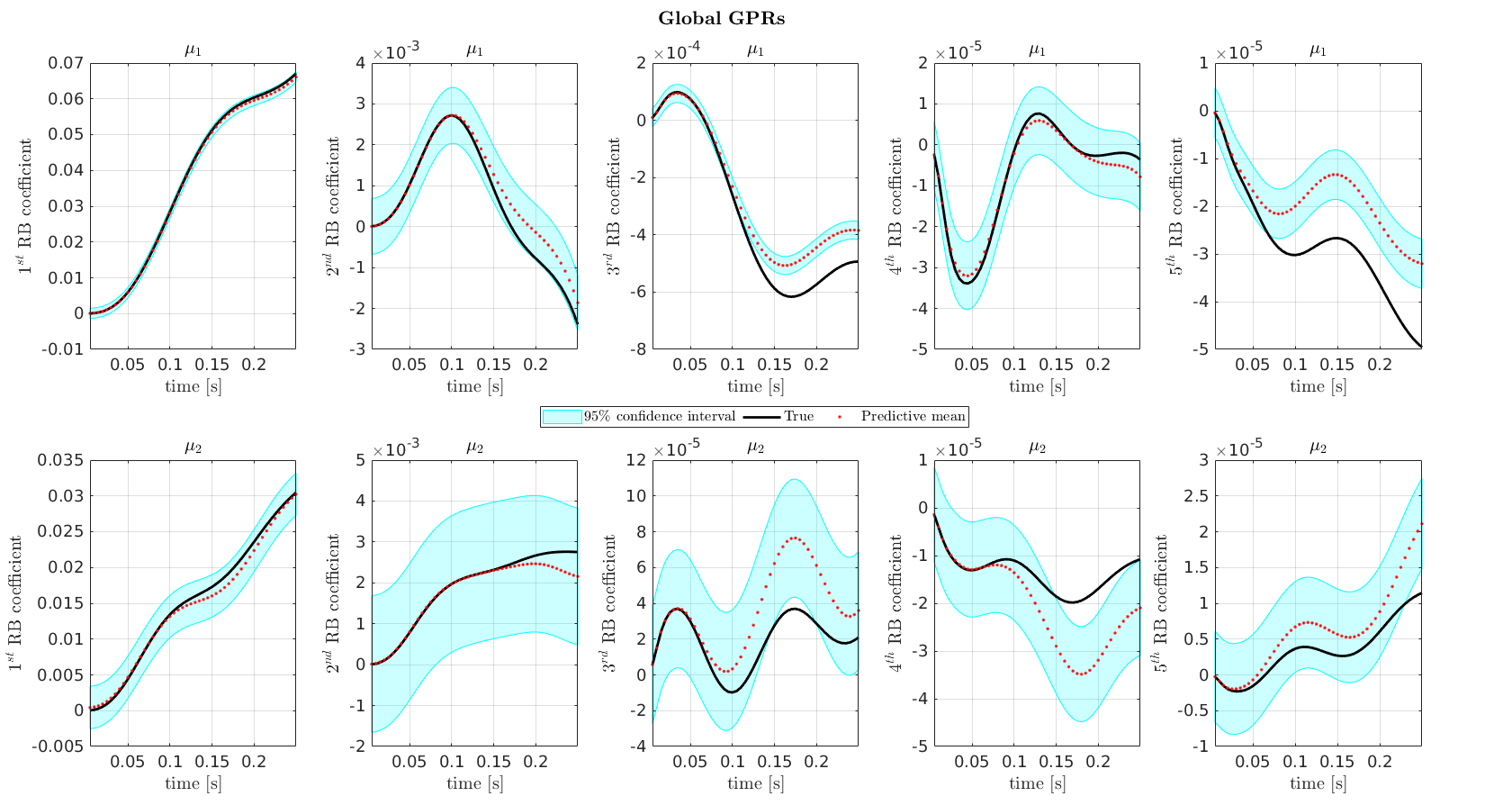}
    \caption{Deformation of a beam. Evolution over time of the exact RB coefficients (black) -- from left to right -- and the corresponding global POD-GPR ROM means (dotted red) for different testing parameters, from top to bottom. Moreover, we report the $95\%$ confidence levels.}
    \label{fig:beam_sampling_mono}
\end{figure}

\begin{figure}[h!]
    \centering
    \includegraphics[width=0.95\textwidth]{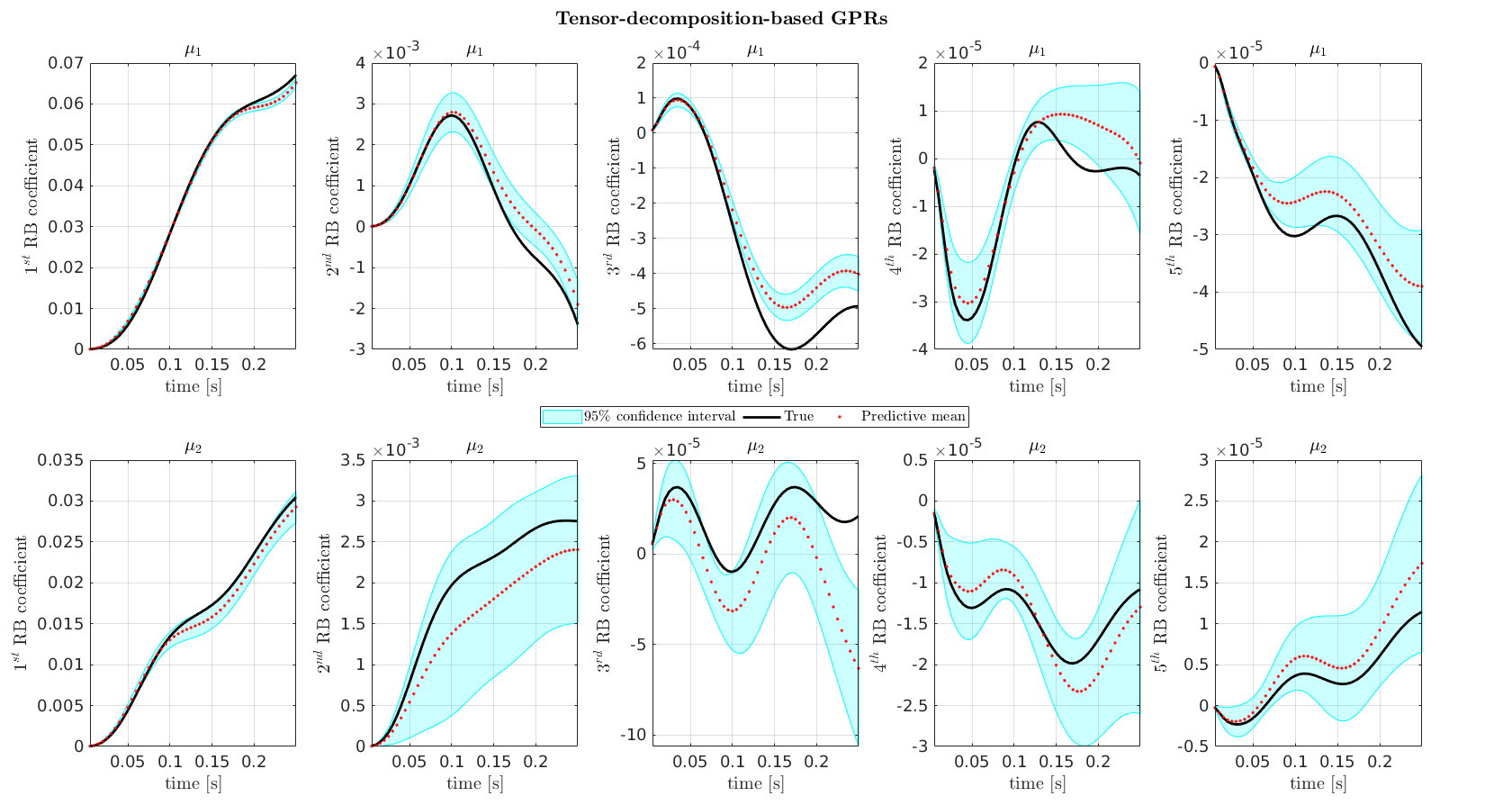}
    \caption{Deformation of a beam. Evolution over time of the exact RB coefficients (black) -- from left to right -- and the corresponding tensor-decomposition-based POD-GPR ROM means (dotted red) for different testing parameters, from top to bottom. Moreover, we report the $95\%$ confidence levels.}
    \label{fig:beam_sampling_segr}
\end{figure}

In Table~\ref{tab:beam_RBcoeff_mono_segr} we report the training time and the RSEs computed when $N_s=50$. The true projection coefficients, the predictive means and their 95\% confidence level computed for two different input parameters are shown in Figure~\ref{fig:beam_sampling_mono} and \ref{fig:beam_sampling_segr} for the global and the tensor-decomposition-based GPR models, respectively. We observe that, in both cases, the errors between the POD-GPR ROMs approximation and the true solution are essentially bounded by the $95\%$ confidence levels.

\begin{table}[h!]
    \centering
    \begin{tabular}{|l|cc|}
        \hline
        & Global GPR & TD-based GPR \\
        \hline
        Training CPU time & $881$ s & $6$ s \\
        $\text{RSE}(q_1)$ & $0.0029$ & $0.0020$ \\
        $\text{RSE}(q_2)$ & $0.0545$ & $0.0535$ \\
        $\text{RSE}(q_3)$ & $0.0732$ & $0.0612$ \\
        $\text{RSE}(q_4)$ & $0.1061$ & $0.1532$ \\
        $\text{RSE}(q_5)$ & $0.1206$ & $0.0705$ \\
        \hline
    \end{tabular}
    \caption{Deformation of a beam. Computational data related to the global and the tensor-decomposition-based regression approaches with ARD-RBF kernels, for $N=5$, $N_s=50$ and $N_t=50$.}
    \label{tab:beam_RBcoeff_mono_segr}
\end{table}

\paragraph{Accuracy of ROM w.r.t FOM.}
The reduced-order solutions $\mathbf{V}\widehat{\mathbf{q}}(t^n;\cdot)$, for $n=1,\dots,N_t$, are reconstructed online to validate the overall accuracy of the regression models with respect to the FOM. First of all, we compute the errors between the high-fidelity solution $\mathbf{u}_h$ and its projection onto the RB space, i.e., $\mathbf{VV}^T\mathbf{u}_h$, averaged over the whole test set, obtaining $\text{tAE}\approx1.27\cdot10^{-5}$ and $\text{tRE}\approx2.01\cdot10^{-3}$. These represent the best (linear) approximation errors in $Col(\mathbf{V})$, and thus the comparative values for our ROMs. From the results reported in Table~\ref{tab:beam_displacement_error}, we observe good accuracy of the GPR models with respect to either the high-fidelity solutions as well as to the RB projection approximations.  We recall that the tensor-decomposition-based approach beats the global one in terms of CPU time required during the offline training. For this test case the FOM requires in average $48$ s to compute the solution dynamics given the parameter vector, whereas the global and the tensor-decomposition-based POD-GPR ROMs are able to reduce the online CPU time up to $364$ and $270$ times, respectively. Furthermore, we compare the accuracy and efficiency of standard POD-Galerkin ROMs equipped with the discrete empirical interpolation method (DEIM) \cite{chaturantabut2010nonlinear} hyper-reduction technique. In particular, we take into account $|\mathcal{I}|= 83$ and $|\mathcal{I}|= 119$ DEIM interpolation points, corresponding to a POD tolerance equal to $10^{-4}$ and $10^{-5}$ on residual snapshots, respectively. The results are reported in Table~\ref{tab:beam_displacement_error}.

\begin{table}[h!]
    \centering
    \begin{tabular}{|l|cccc|}
        \hline
        & \multicolumn{2}{c}{POD-Galerkin-DEIM} & Global GPR & TD-based GPR \\
        & $|\mathcal{I}|= 83$ & $|\mathcal{I}|= 119$ & & \\
        \hline
        Online CPU time & $11.4$ s & $13.4$ s & $0.13$ s & $0.18$ s \\
        Speed-up & $4.2$ & $3.6$ & $364$ & $270$ \\
        $\text{mean}_{\boldsymbol{\mu}}~\text{tAE}(\boldsymbol{\mu})$ & $2.3\cdot10^{-4}$ & $1.2\cdot10^{-4}$ & $4.4\cdot10^{-4}$ & $4.1\cdot10^{-4}$ \\
        $\text{mean}_{\boldsymbol{\mu}}~\text{tRE}(\boldsymbol{\mu})$ & $1.6\cdot10^{-2}$ & $1.3\cdot10^{-2}$ & $1.0\cdot10^{-1}$ & $2.1\cdot10^{-2}$ \\
        \hline
    \end{tabular}
    \caption{Deformation of a beam. Efficiency and accuracy of POD-Galekin-DEIM ROM, global POD-GPR ROM and tensor-decomposition-based POD-GPR ROM.}
    \label{tab:beam_displacement_error}
\end{table}

\begin{figure}[h!]
    \centering
    \includegraphics[width=0.95\textwidth]{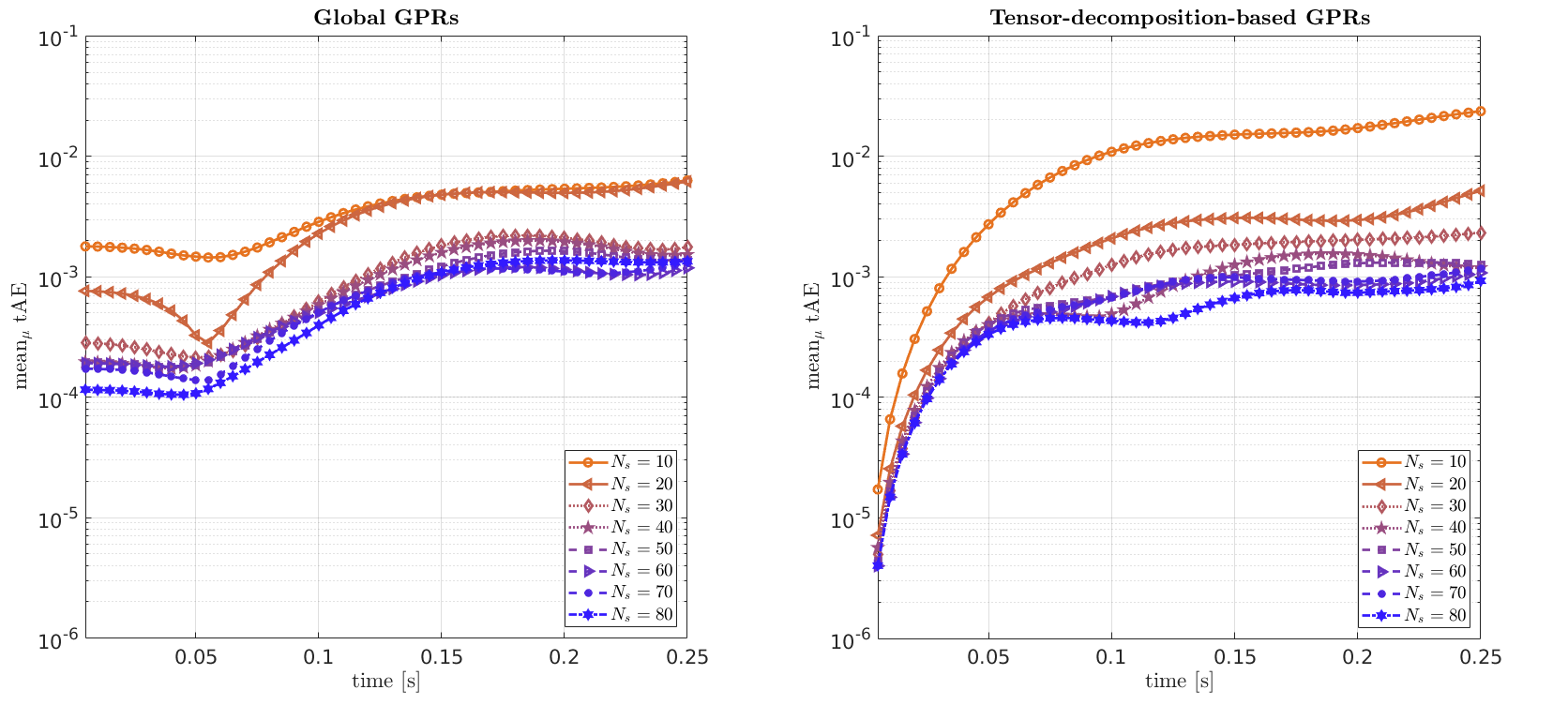}
    \caption{Deformation of a beam. Average absolute error of the global (left) and tensor-de\-com\-po\-si\-tion-based (right) POD-GPR ROMs approximations over the size of the training set.}
    \label{fig:beam_displacement_tAE}
\end{figure}

\begin{figure}[h!]
    \centering
    \includegraphics[width=0.95\textwidth]{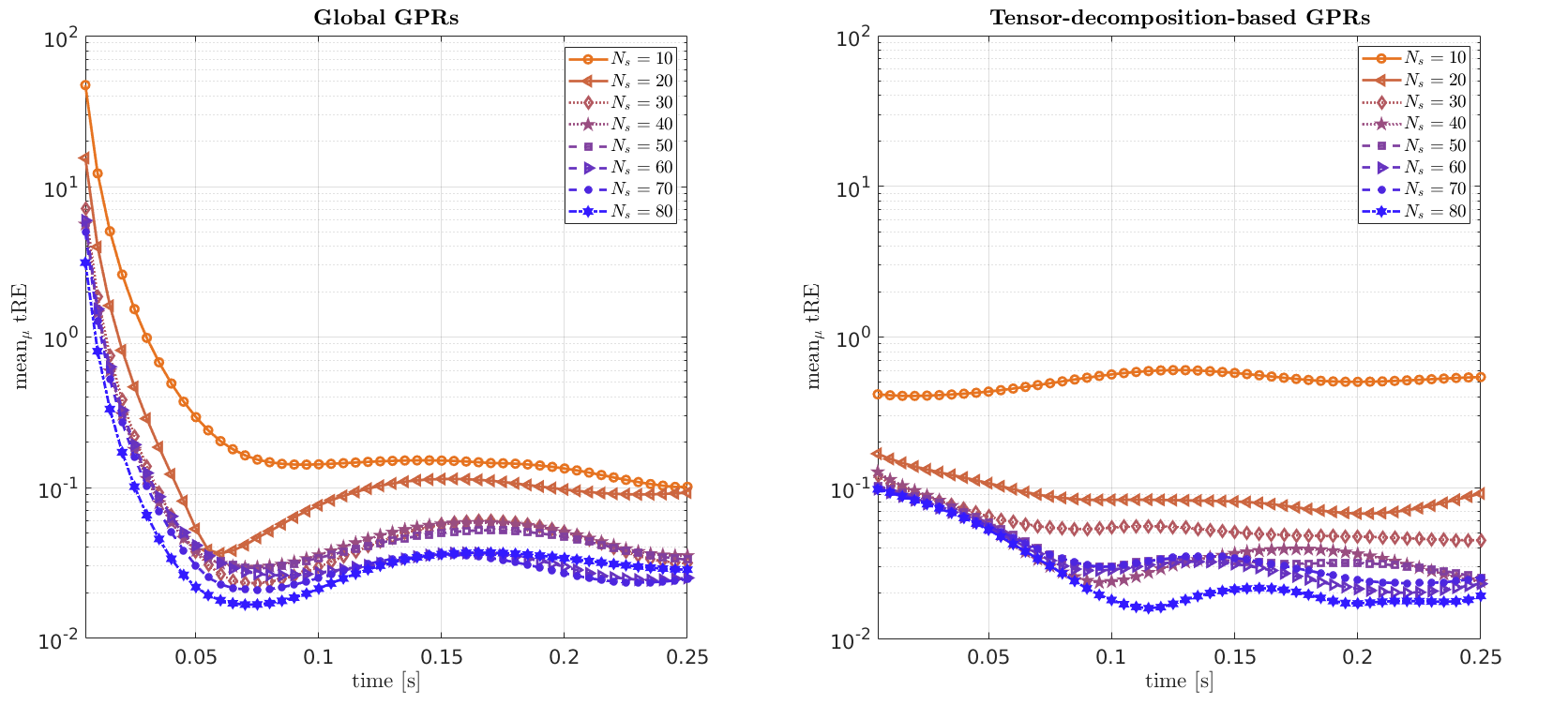}
     \vspace{-2mm}
    \caption{Deformation of a beam. Average relative error of the global (left) and tensor-de\-com\-po\-si\-tion-based (right) POD-GPR ROMs approximations over time, varying the size of the training set.}
    \label{fig:beam_displacement_tRE}
\end{figure}

\begin{figure}[h!]
\vspace{-2mm}
    \centering
    \includegraphics[width=0.95\textwidth]{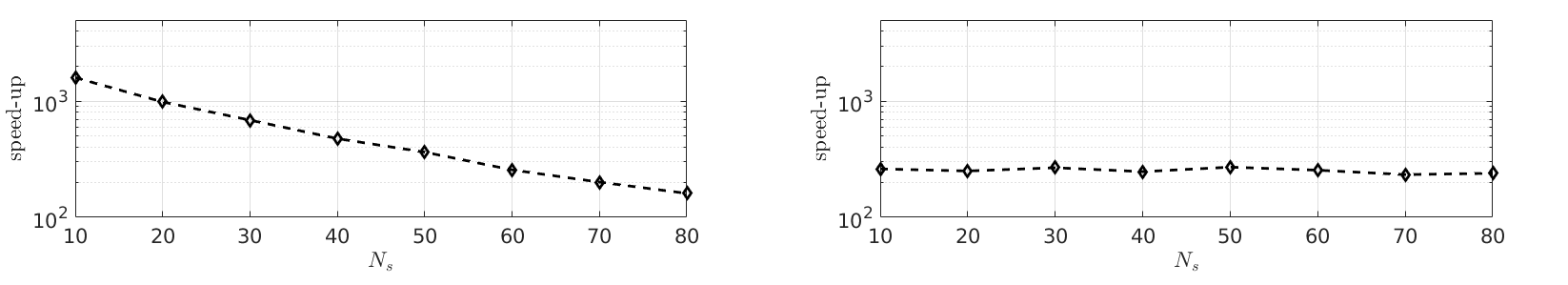}
    \vspace{-2mm}
    \caption{Deformation of a beam. Average computational speed-up of the global (left) and tensor-de\-com\-po\-si\-tion-based (right) POD-GPR ROMs approximations over time, varying the size of the training set, with respect to the full-order solution.}
    \label{fig:beam_displacement_time}
\end{figure}

\begin{figure}[h!]
    \includegraphics[width=\textwidth]{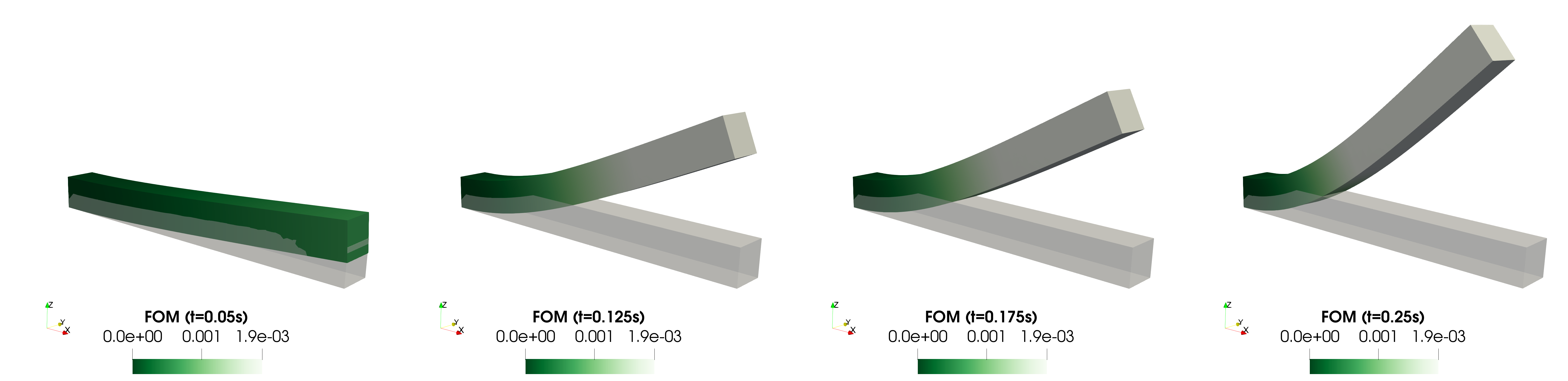}\\
    \includegraphics[width=\textwidth]{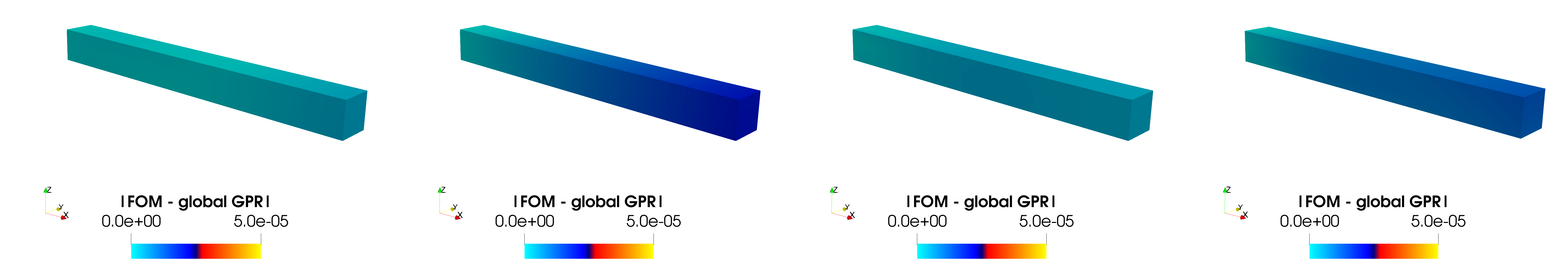}\\
    \includegraphics[width=\textwidth]{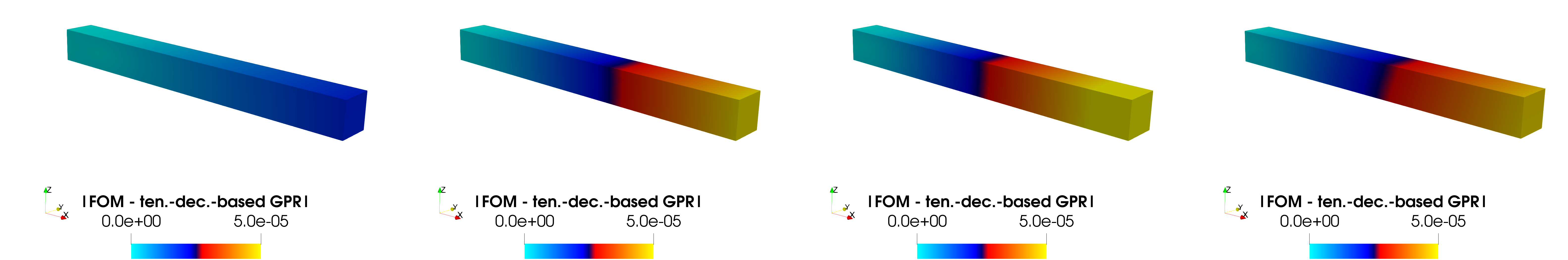}
    \vspace{-3mm}
    \caption{Deformation of a beam. FOM displacement (top), and approximation error of the global (middle) and tensor-decomposition-based (bottom) POD-GPR ROM, at different time instances, for $\boldsymbol{\mu} = [4.20, 2.63, 2.51, 4.54, 4.46, 2.81, 27.25~\text{kPa}, 2.51~\text{kPa}, 0.0052~\text{kPa}]$.}
    \label{fig:beam_mu15}
\end{figure}

Figures~\ref{fig:beam_displacement_tAE} and \ref{fig:beam_displacement_tRE} show the evolution in time of the mean (over the test set) absolute and relative errors, respectively, computed with respect the FOM solution for different number of training samples $N_s$. Moreover, the corresponding computational speed-ups are reported in Figure~\ref{fig:beam_displacement_time}. We observe that both approaches perform well and with comparable accuracy when the training dataset is sufficiently large. In particular, global GPRs are able to achieve higher computational speed-ups than tensor-decomposition-based ones, although they strictly depend on $N_s$. The FOM displacement and the approximation error of the POD-GPR ROMs at different time instances are reported in Figure~\ref{fig:beam_mu15}, for a given values of the input vector.

\vspace{3mm}

\begin{figure}[h!]
    \centering
    \includegraphics[width=0.95\textwidth]{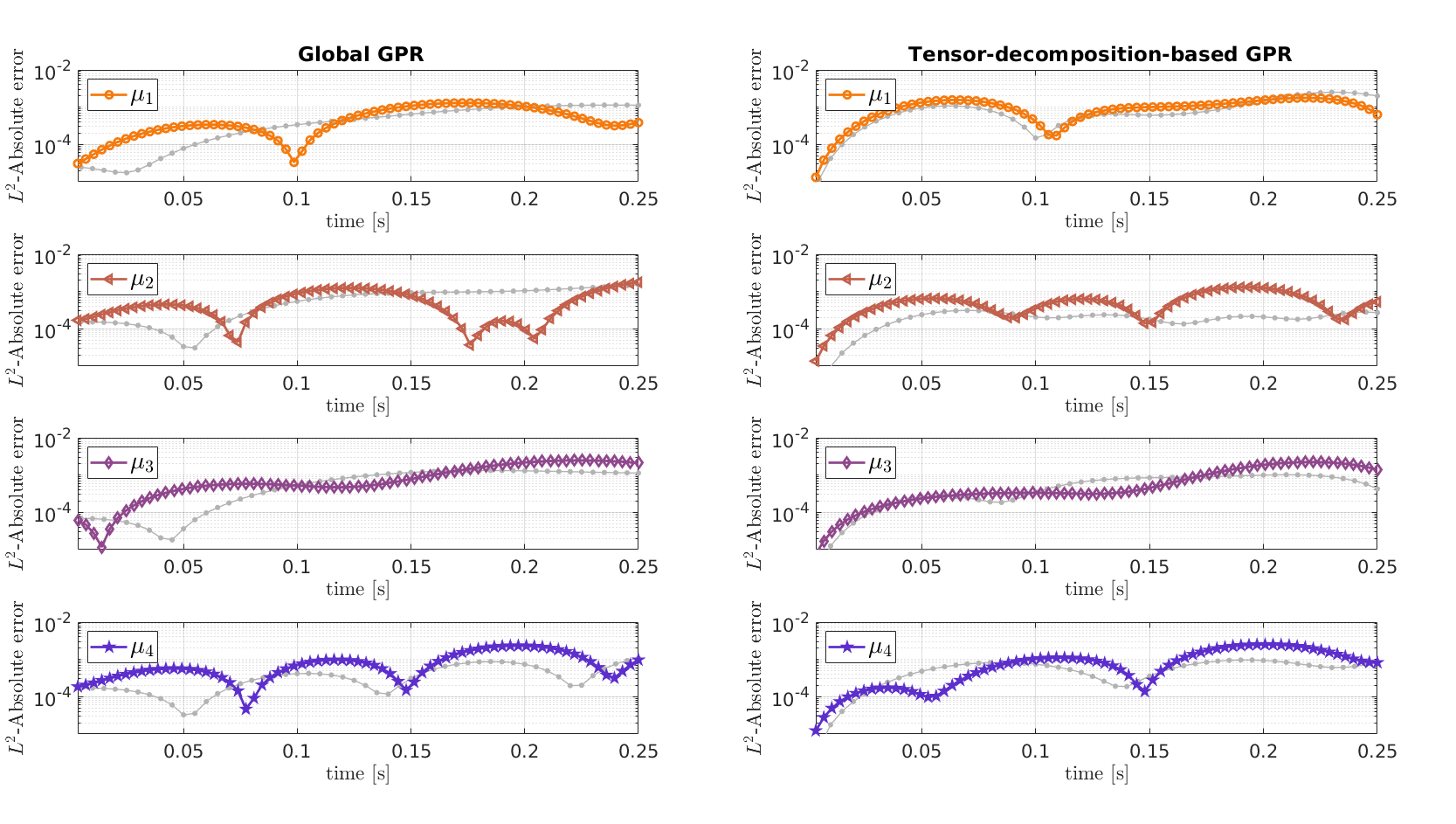}
    \caption{Deformation of a beam. Evolution of the time-average $L^2$-absolute error of the global (left) and tensor-decomposition-based (right) POD-GPR ROMs, using $\Delta t = 3.5\cdot10^{-3}$~s.}\label{fig:beam_displacement_error_time_generalization}
\end{figure}

So far we have always considered during the online stage the same time instances used for training the GPRs, since out main focus is on the influence of the parameters on the deformation of the body. However, GPRs can generalize to unseen values also for the time component, as shown by the computational results reported in Table~\ref{tab:beam_displacement_error_time_generalization} when a time step size $\Delta t = 3.5\cdot10^{-3}$~s is used online. In this case the total number of time stops is $N_t^{test}=71$, and the FOM requires $75$~s in average for a given parameter vector. The corresponding time-average $L^2$-absolute error is reported in Figure~\ref{fig:beam_displacement_error_time_generalization} for different values of the input parameter. Within each plot we also report, in light gray color, the error computed using the same time step used during training, showing good prediction properties of the POD-GPR ROMs with respect both time and parameters.

\begin{table}[h!]
    \centering
    \begin{tabular}{|l|cc|}
        \hline
        & Global GPR & TD-based GPR \\
        \hline
        Online CPU time & $0.18$ s & $0.20$ s \\
        Speed-up & $385$ & $366$ \\
        $\underset{\boldsymbol{\mu}}{\text{mean}}~tAE$ & $1.1\cdot10^{-3}$ & $1.0\cdot10^{-3}$ \\
        $\underset{\boldsymbol{\mu}}{\text{mean}}~tRE$ & $3.9\cdot10^{-1}$ & $7\cdot10^{-2}$ \\
        \hline
    \end{tabular}
    \caption{Deformation of a beam. Efficiency and accuracy of global and tensor-decomposition-based POD-GPR ROMs, computed over $20$ testing parameters, using $\Delta t = 3.5\cdot10^{-3}$~s.}\label{tab:beam_displacement_error_time_generalization}
\end{table}	

\subsubsection{Test 2: active contraction of a truncated ellipsoid} \label{sec:prolate_results}

In this second test case, a larger number of RB functions is required to achieve the same POD accuracy with respect to the previous problem (see Table~\ref{tab:prolate_POD}), given the slower decay of the singular values, thus increasing the complexity of the surrogate model due to the higher number of GPs that must be trained.

\begin{table}[h!]
    \centering
    \begin{tabular}{|c|ccccc|}
        \hline
        $\boldsymbol{\varepsilon}_{\bf POD}$ & $10^{-3}$ & $5\cdot10^{-4}$ & $10^{-4}$ & $5\cdot10^{-5}$ & $10^{-5}$\\
        $N$ & 20 & 30 & 71 & 100 & 202 \\
        \hline
    \end{tabular}
    \caption{Active contraction of a truncated ellipsoid. Dimension of RB subspaces  for different POD tolerances.}
    \label{tab:prolate_POD}
\end{table}

\paragraph{Predictions of GPRs on RB coefficients.}
Once fixed the RB dimension equal to $N=30$, we analyze the fidelity of approximations computed using the the global and the tensor-decomposition-based POD-GPR ROMs built using different covariance functions, trained on $N_s=50$ parameter samples. Regarding the tensor-decomposition-based approach, time- and parameter-modes for each RB coefficient
are obtained by computing the truncated SVD, obtaining $3\leq N_\ell\leq9$. Based on the observation of Section~\ref{sec:beam_results}, we only rely on the standardization scaling of the training and test sets.

\begin{figure}[h!]
    \centering
    \includegraphics[width=0.95\textwidth]{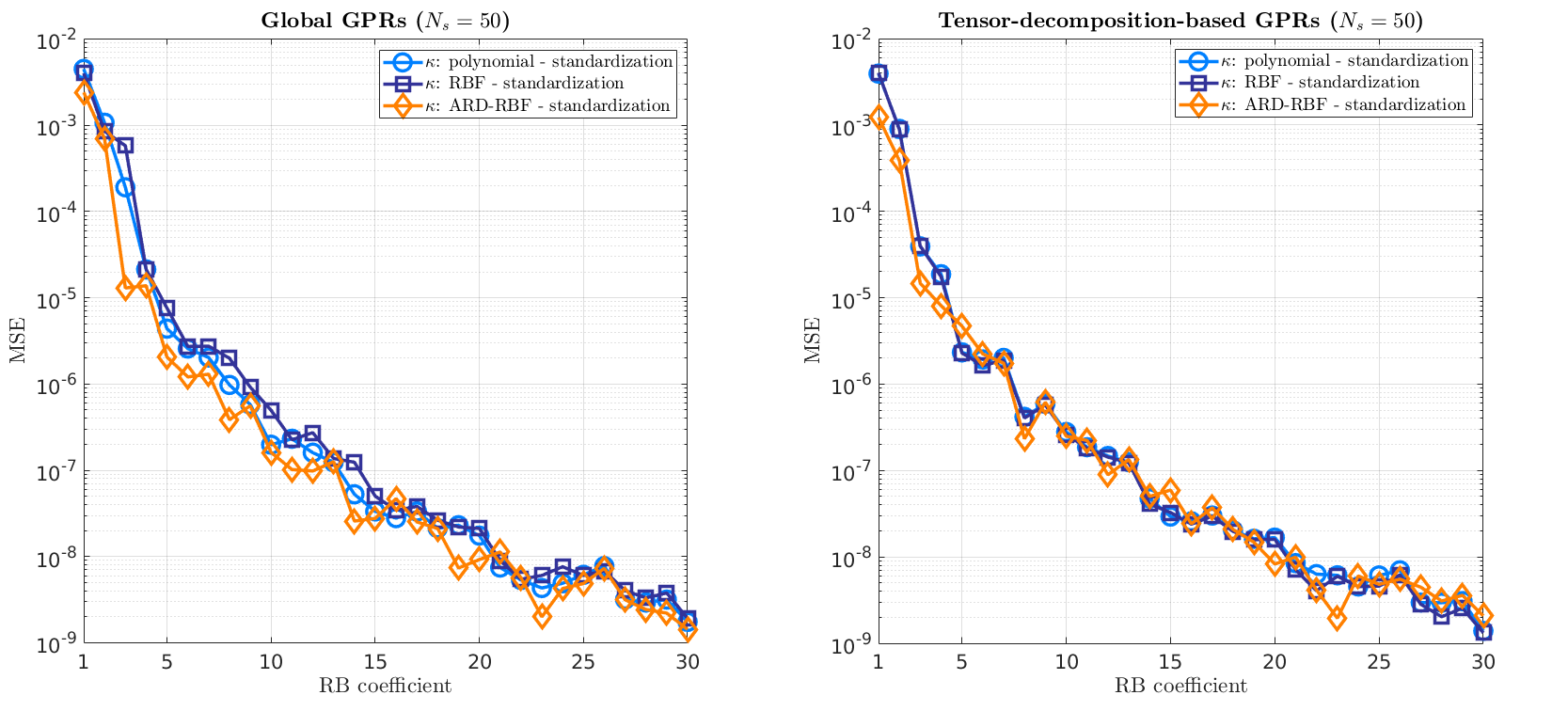}
    \caption{Active contraction of a truncated ellipsoid. Mean square error of global (left) and tensor-decomposition-based (right) POD-GPR ROMs computed using different kernel types, for $N=30$.}
    \label{fig:prolate_Ns_kernel_scaling}
\end{figure}

From the MSE reported in Figure~\ref{fig:prolate_Ns_kernel_scaling}, we observe that the choice of the kernel function has a smaller impact on the accuracy of the GPR predictions with respect to the previous test. However, since the ARD-RBF kernel shows better agreement to the true values of the first RB coefficient than the other kernels, we restrict ourselves to this choice for the GPs covariance function, unless otherwise specified.

\begin{figure}[h!]
    \centering
    \includegraphics[width=0.95\textwidth]{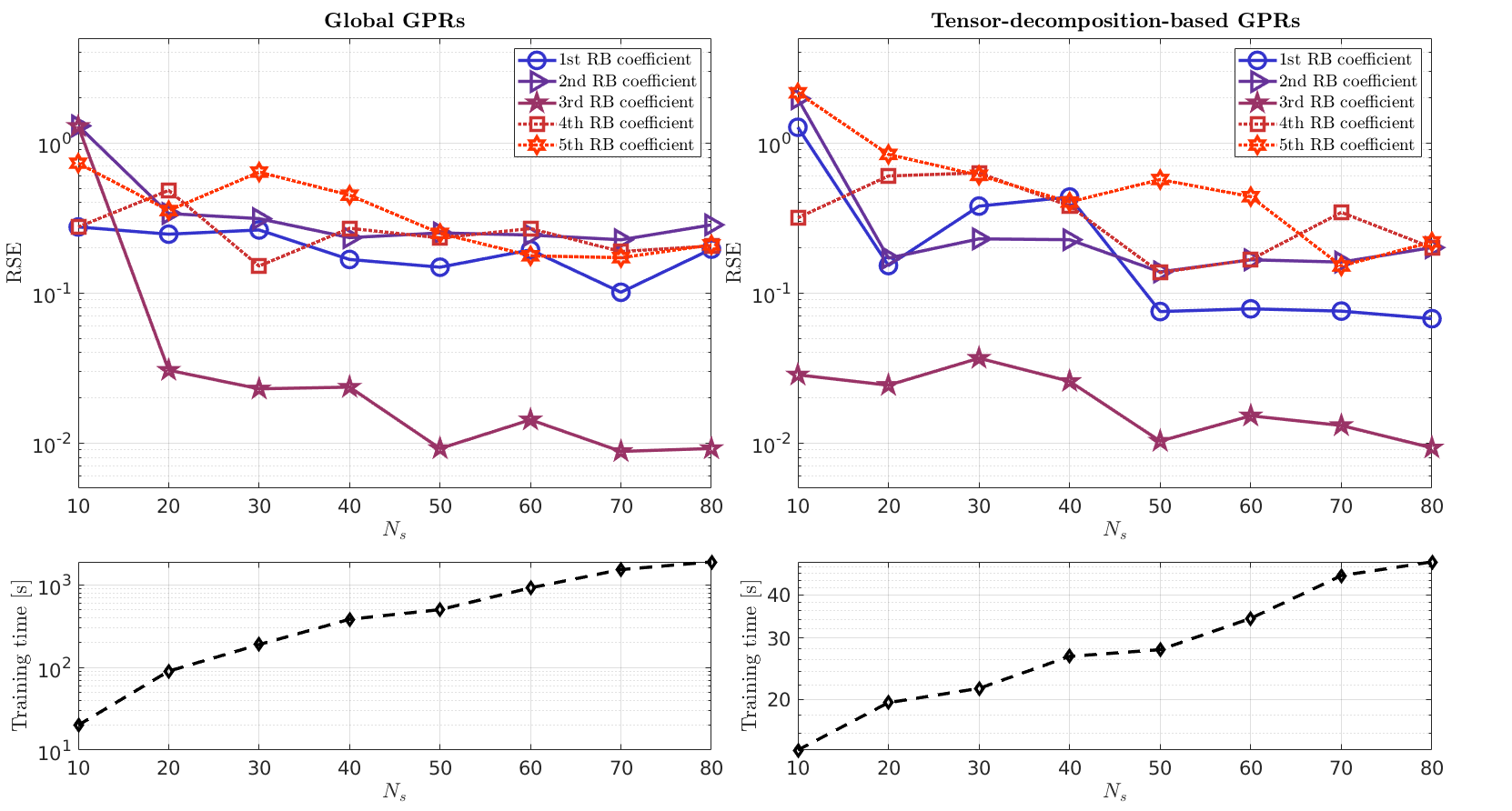}
    \caption{Active contraction of a truncated ellipsoid. Relative squared error -- associated with the first $5$ RB coefficients -- of the global (top left) and the tensor-decomposition-based (top right) POD-GPR ROMs over the size of the training set. Corresponding training times (bottom) are also reported. These results have been obtained taking into account a testing set of $N_s^{test}=20$ parameters, unseen during training.}
    \label{fig:prolate_convergence_mono_segr}
\end{figure}

Hence, we analyze the accuracy and the efficiency of the POD-GPR ROMs by increasing the size of the training set from $N_s=10$ up to $80$. We point out that the global approach requires approximately $4$~h to train the GPs on $N_s=50$ sample parameters, probably due to the high FOM dimension $N_h$, in addition to the larg number of RB coefficients. To reduce this burden without restoring to parallel computing, from now on we restrict the size of the training set by keeping only the data associated with time instant $t^{5i}$, for $i=1,\dots,\lfloor N_t/5 \rfloor$. From Figure~\ref{fig:prolate_convergence_mono_segr} we observe that a dataset compromising $N_s=40$ or $N_s=50$ sampled parameters represents a good trade-off between accuracy and CPU training times. 

\begin{figure}[h!]
    \centering
    \includegraphics[width=0.95\textwidth]{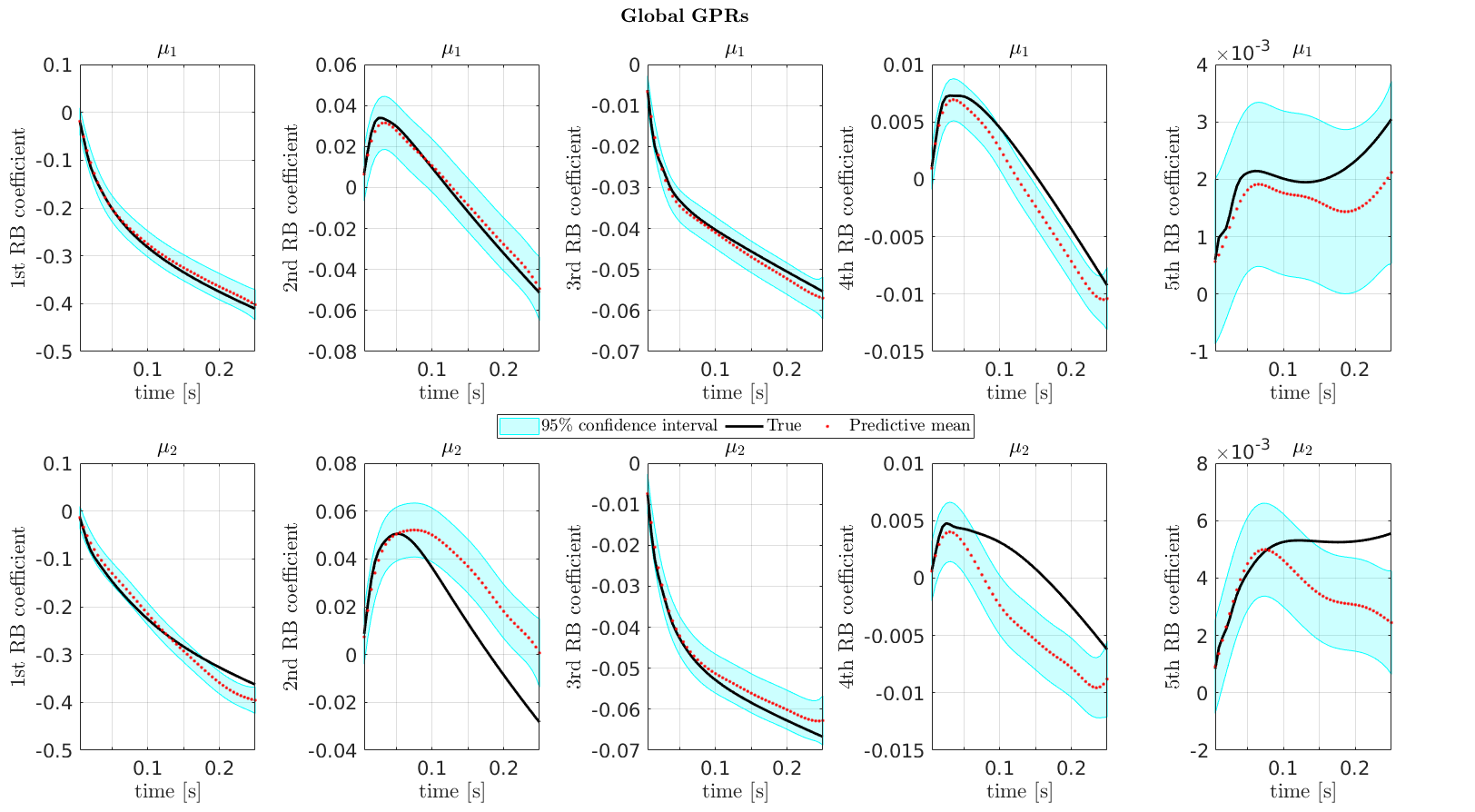}
    \caption{Active contraction of a truncated ellipsoid. Evolution over time of the exact RB coefficients (black) -- from left to right -- and the corresponding global POD-GPR ROM means (dotted red) for different testing parameters, from top to bottom. Moreover, we report the $95\%$ confidence levels.}
    \label{fig:prolate_sampling_mono}
\end{figure}

\begin{figure}[h!]
    \centering
    \includegraphics[width=0.95\textwidth]{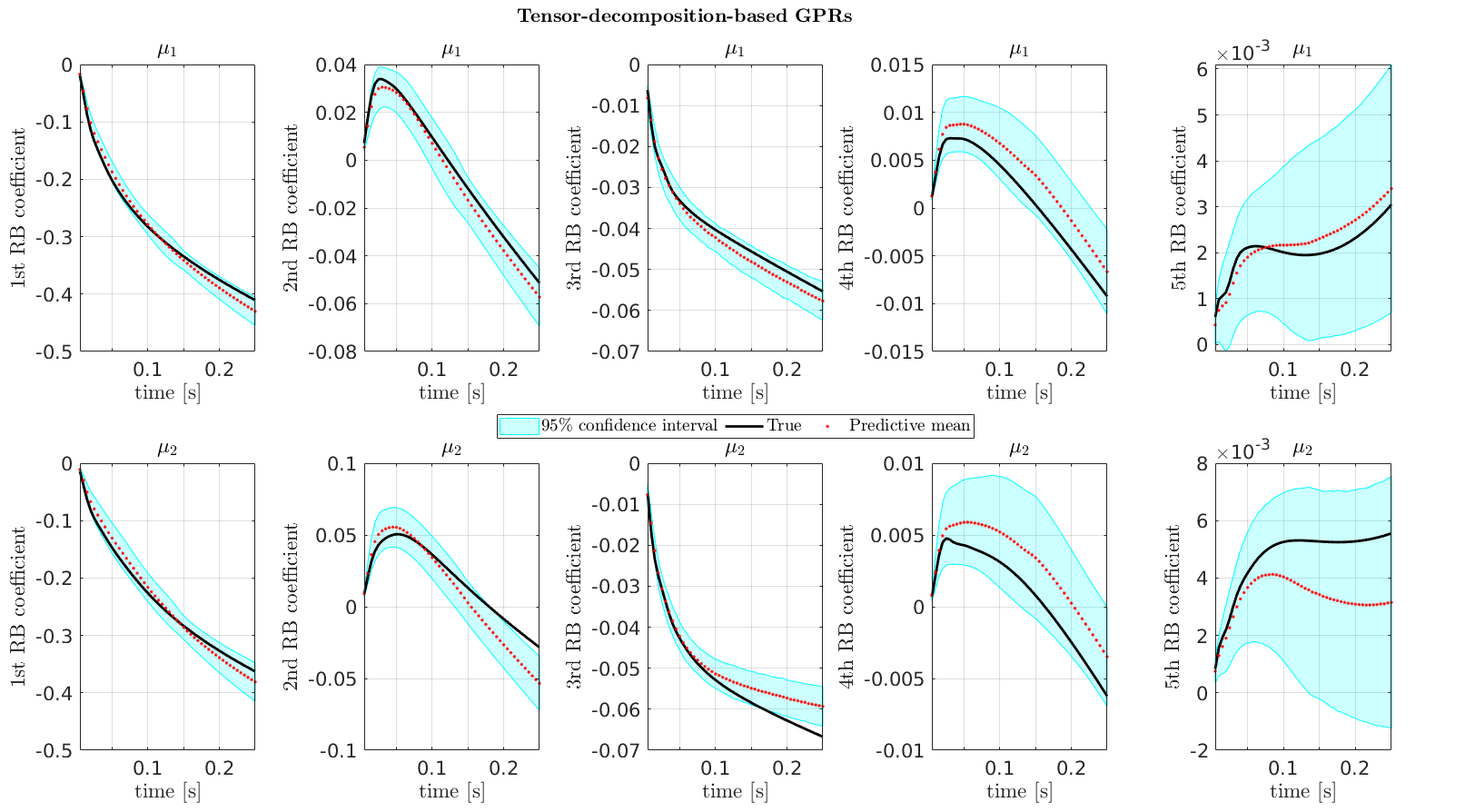}
    \caption{Active contraction of a truncated ellipsoid. Evolution over time of the exact RB coefficients (black) -- from left to right -- and the corresponding tensor-decomposition-based POD-GPR ROM means (dotted red) for different testing parameters, from top to bottom. Moreover, we report the $95\%$ confidence levels.}
    \label{fig:prolate_sampling_segr}
\end{figure}

For the RB coefficients $q_1,\dots,q_5$ we report in Table~\ref{tab:prolate_RBcoeff_mono_segr} the computational data obtained with $N_s=50$, whilst the corresponding true projection coefficients, predictive means and 95\% confidence level are depicted in Figures~\ref{fig:prolate_sampling_mono} and \ref{fig:prolate_sampling_segr}. The two approaches are almost comparable in terms of accuracy, whilst the tensor-decomposition-based approach shows faster training times.
    
\begin{table}[h!]
    \centering
    \begin{tabular}{|l|cc|}
        \hline
        & Global GPR & TD-based GPR \\
        \hline
        Training CPU time & $501$ s & $28$ s \\
        $\text{RSE}(q_1)$ & $0.1487$ & $0.0753$ \\
        $\text{RSE}(q_2)$ & $0.2517$ & $0.1384$ \\
        $\text{RSE}(q_3)$ & $0.0092$ & $0.0103$ \\
        $\text{RSE}(q_4)$ & $0.2331$ & $0.1369$ \\
        $\text{RSE}(q_5)$ & $0.2501$ & $0.5711$ \\
        \hline
    \end{tabular}
    \caption{Active contraction of a truncated ellipsoid. Computational data related to the global and the tensor-decomposition-based regression approaches with ARD-RBF kernels, for $N=30$, $N_s=50$ and $N_t=50$.}
    \label{tab:prolate_RBcoeff_mono_segr}
\end{table}	

\begin{figure}[h!]
    \centering
    \includegraphics[width=0.95\textwidth]{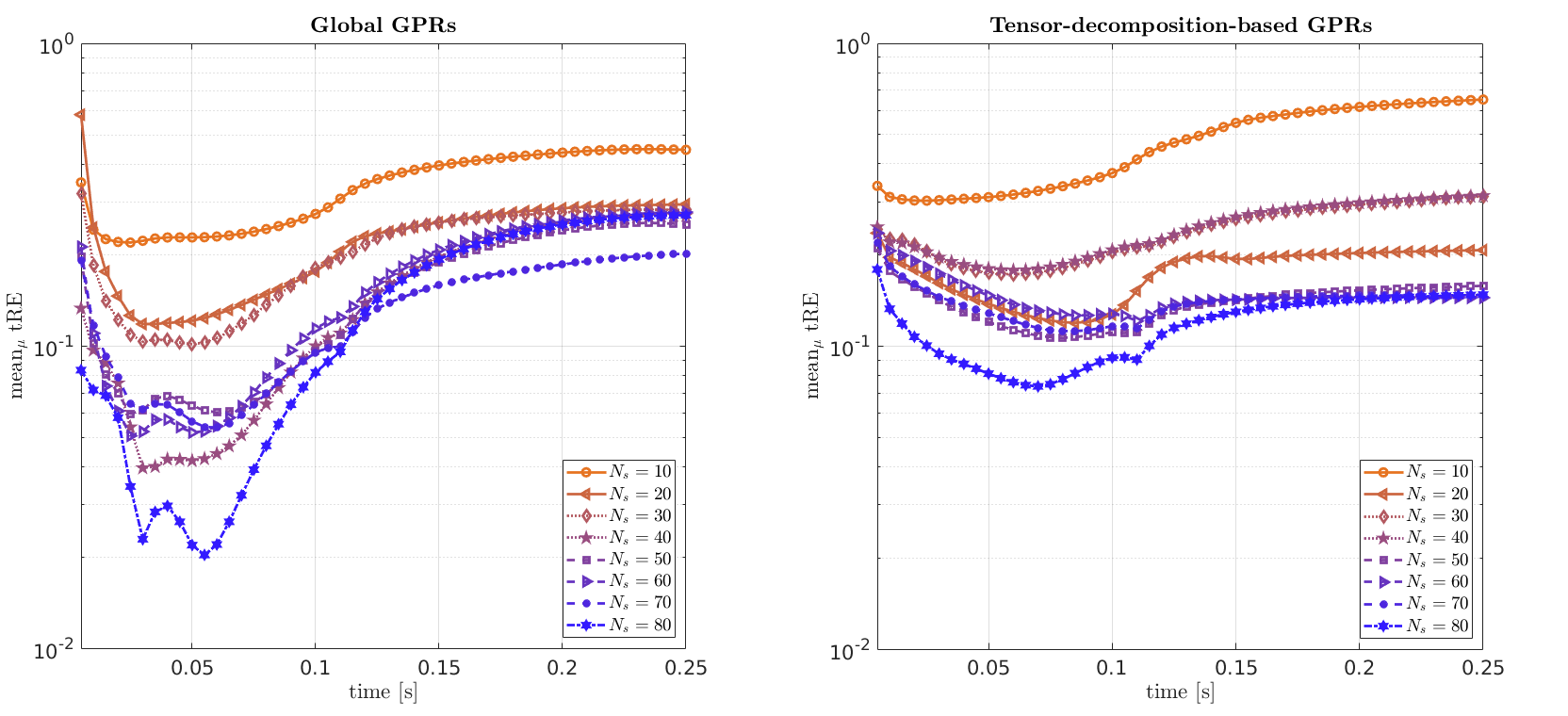}
    \caption{Active contraction of a truncated ellipsoid. Average relative error of the global (left) and tensor-de\-com\-po\-si\-tion-based (right) POD-GPR ROMs approximations over time, varying the size of the training set.}
    \label{fig:prolate_displacement_tRE}
\end{figure}

\begin{figure}[h!]
    \centering
    \includegraphics[width=0.95\textwidth]{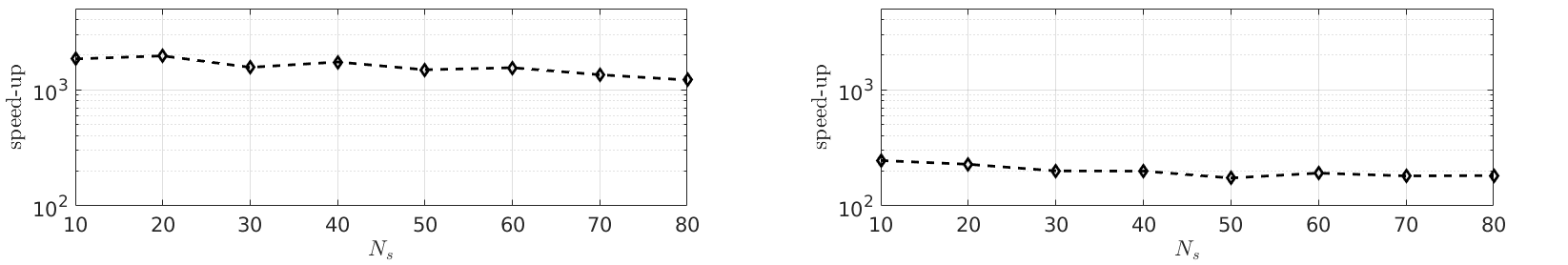}
    \caption{Active contraction of a truncated ellipsoid. Average computational speed-up of the global (left) and tensor-de\-com\-po\-si\-tion-based (right) POD-GPR ROMs approximations over the size of the training set, with respect to the full-order solution.}
    \label{fig:prolate_displacement_time}
\end{figure}

\begin{figure}[h!]
    \includegraphics[width=\textwidth]{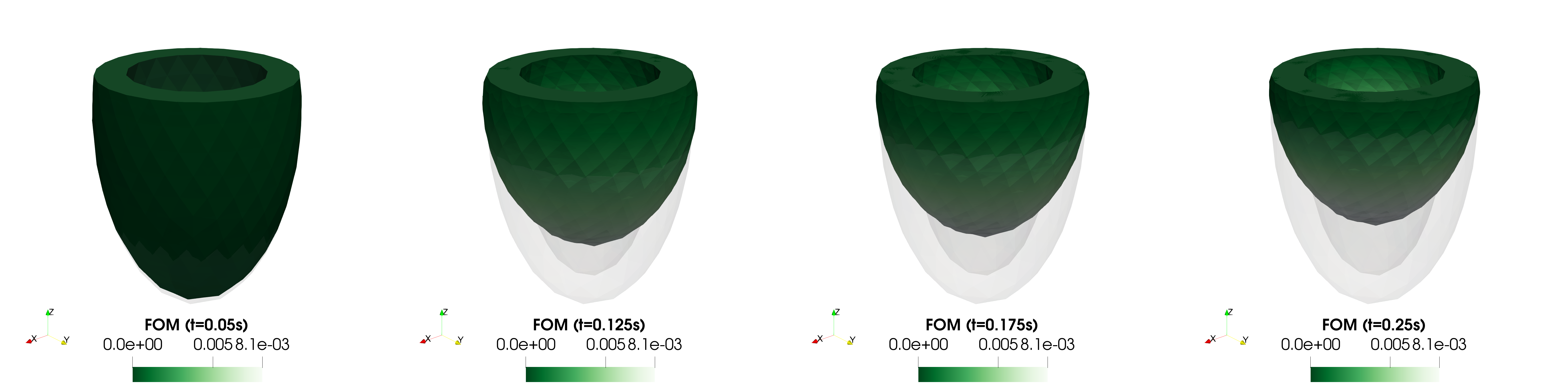}\\
    \includegraphics[width=\textwidth]{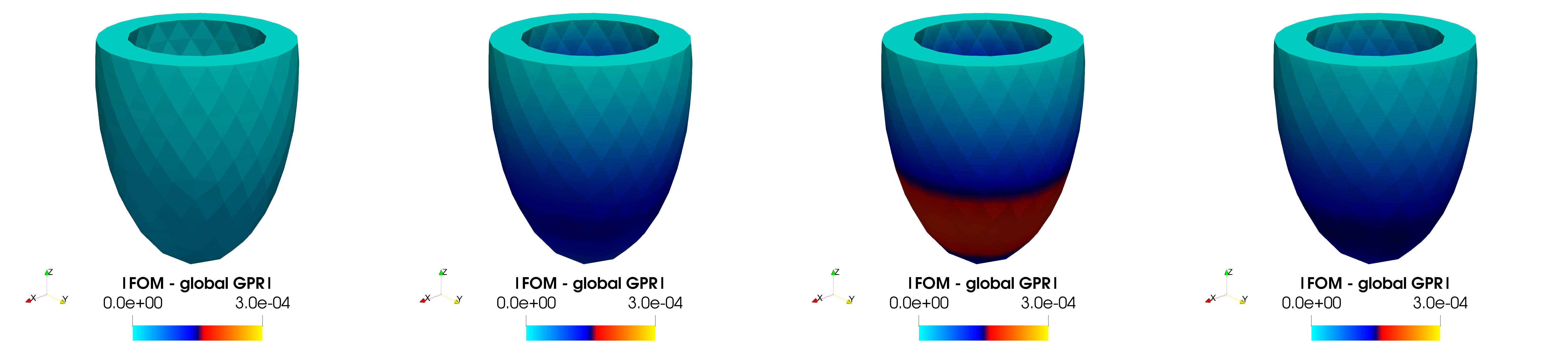}\\
    \includegraphics[width=\textwidth]{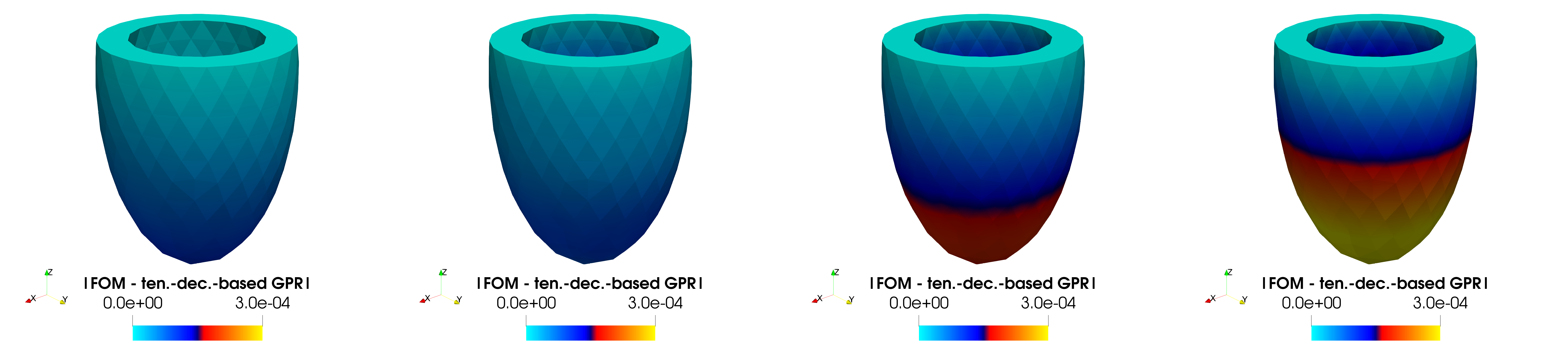}
    \vspace{-3mm}
    \caption{Active contraction of a truncated ellipsoid. FOM displacement (top), and approximation error of the global (middle) and tensor-decomposition-based (bottom) POD-GPR ROM, at different time instances, for $\boldsymbol{\mu} = [7.85, 1.81, 2.21, 4.11, 4.29, 1.76, 51.30~\text{kPa}, 1.92~\text{kPa}, 64.94~\text{kPa}, 14.73~\text{kPa}, -98.53^\circ, 104.72^\circ]$.}
    \label{fig:prolate_mu1}
\end{figure}

\paragraph{Accuracy of ROM w.r.t. FOM.}
Finally, we compare the accuracy of the surrogate GPR models on the approximation of the whole high-fidelity displacement field, that is how $\mathbf{V}\hat{\mathbf{q}}(\cdot;\boldsymbol{\mu})$ differs from $\mathbf{u}_h(\cdot;\boldsymbol{\mu})$, as well as its computational efficiency. We recall that the FOM requires $6$ min $11$ s in average to compute the solution dynamics for a given parameter. The results obtained are reported in Figure~\ref{fig:prolate_displacement_tRE} and \ref{fig:prolate_displacement_time}, and summarized in Table~\ref{tab:prolate_displacement_error}. As already observed, the speed-up achieved online by the global POD-GPR ROMs is the highest between the two approaches. Figure~\ref{fig:prolate_mu1} shows the FOM displacement and the approximation error of the POD-GPR ROMs at different time instances, for a given values of the input vector.

\begin{table}[h!]
    \centering
    \begin{tabular}{|l|cc|}
        \hline
        & Global GPR & TD-based GPR \\
        \hline
        Online CPU time & $0.25$ s & $2.1$ s \\
        Speed-up & $1487$ & $174$ \\
        $\text{mean}_{\boldsymbol{\mu}}~\text{tAE}(\boldsymbol{\mu})$ & $3.1\cdot10^{-2}$ & $2.8\cdot10^{-2}$ \\
        $\text{mean}_{\boldsymbol{\mu}}~\text{tRE}(\boldsymbol{\mu})$ & $1.6\cdot10^{-1}$ & $1.4\cdot10^{-1}$ \\
        \hline
    \end{tabular}
    \caption{Active contraction of a truncated ellipsoid. Efficiency and accuracy of global POD-GPR ROM and tensor-decomposition-based POD-GPR ROM.}
    \label{tab:prolate_displacement_error}
\end{table}

\section{Multi-query problems in UQ using POD-GPR ROMs} \label{sec:UQ}

Mathematical models describing physical phenomena depend on a (large) set of parameters that characterize different physical responses, and whose knowledge is severely limited due difficulties in performing measurements. However, only a smaller subset of inputs have a great impact on specific outputs quantities, that are scalar-valued functions of the displacement field, for the application at hand. Hence, it is important to identify which factors should be properly measured and which others can be arbitrarily fixed without affecting the selected QoIs. Another key task when dealing with parameterized models is to infer about the unknown parameters based on observed data. We show how, using the mean of the posterior distribution of the GPRs to make predictions, we are able to efficiently address the solution to both types of problems, i.e., sensitivity analysis and parameter estimation, which require the repeated evaluation of the input-output map.

\begin{figure}[h!]
    \centering
    \includegraphics[width=\textwidth]{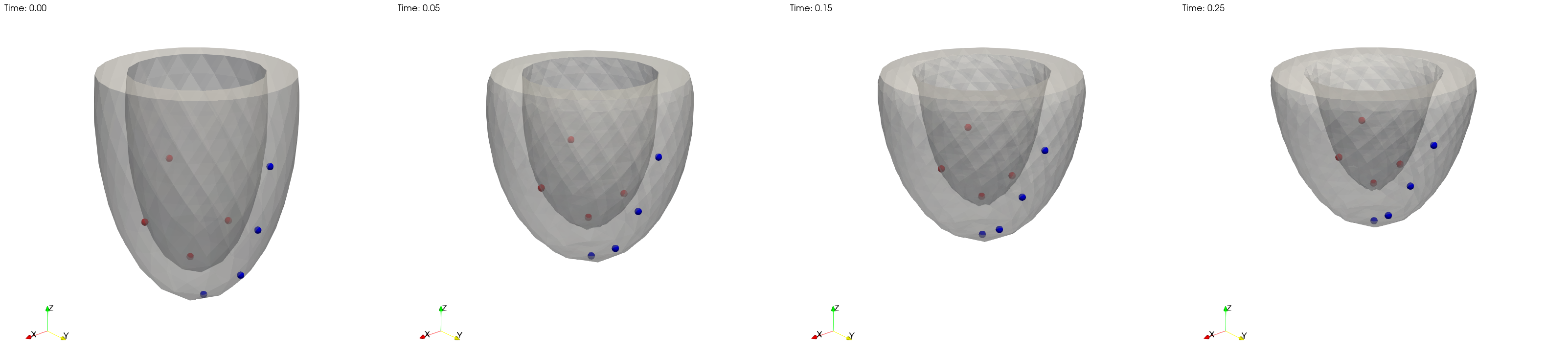}
    \caption{Active contraction of a truncated ellipsoid. Tissue displacement is measured at the epicardium (light red) and at the endocardium (dark red), at three different time instances.}
    \label{fig:prolate_QoI}
\end{figure}

In the following we address the solution to these multi-query problems on the numerical test cases reported in Sections~\ref{sec:beam} and \ref{sec:prolate}. For the deforming beam, as $k$-th QoI we consider the displacement along the $z$-axis of a prescribed point with coordinates $\bar{\boldsymbol{X}} = (10^{-2}\text{m},5\cdot10^{-2}\text{m},5\cdot10^{-2}\text{m})$ in the reference domain $\Omega_0$, measured at time $t=j_k\Delta t$, that is
\begin{equation*}
    y_k(\boldsymbol{\mu})= \mathbf{u}(\bar{\boldsymbol{X}},j_k\Delta t;\boldsymbol{\mu}) \cdot \mathbf{e}_z,
\end{equation*}
for $j_k = 10,30,50$. In the case of the active contraction of a truncated ellipsoid, we take into account the horizontal and vertical displacements of different points located on the outer and inner surfaces and reposted in Figure~\ref{fig:prolate_QoI}, acquired at different time instances, such that the outputs are given by
\begin{align*}
    y_k(\boldsymbol{\mu})&= \mathbf{u}(\bar{\boldsymbol{X}}_{i_k},j_k\Delta t;\boldsymbol{\mu}) \cdot \mathbf{e}_x\\
    y_{k+1}(\boldsymbol{\mu})&= \mathbf{u}(\bar{\boldsymbol{X}}_{i_k},j_k\Delta t;\boldsymbol{\mu}) \cdot \mathbf{e}_y\\
    y_{k+2}(\boldsymbol{\mu})&= \mathbf{u}(\bar{\boldsymbol{X}}_{i_k},j_k\Delta t;\boldsymbol{\mu}) \cdot \mathbf{e}_z
\end{align*}
for $i_k=1,\dots,8$ and $j_k=10,30,50$. With an abuse of notation, we denote as $\boldsymbol{\mu}$ (or $y$) the parameter vector both as random variable (or model output random variable) and as its outcome. To perform the UQ studies we employ the Dakota toolkit \cite{adams2020dakota}, exploiting its Python direct interface.

\subsection{Sensitivity Analysis (SA)}\label{sec:SA}

First of all, we perform global SA, which allows us to rank the parameters according to their influence on the QoIs and thus gain useful information. To begin, we employ the screening method of Morris' elementary effects~\cite{Morris1991}, which provides a qualitative analysis relying only on a small number of numerical simulations. In this case, we compare the results obtained using the GPR models (that, we recall, are trained on a dataset built from $N_s=50$ sampled parameters, and suitable scaled using the standardization technique) with the high-fidelity ones. Hence, we perform a variance-based global SA by computing Sobol' first order and total effects indices~\cite{sobol1990sensitivity}, providing a more quantitative ranking criterion, although a much larger number of input-output evaluations is required. Since multiple outputs are analyzed by conducting independent sensitivity studies, we restrict the formal introduction to a univariate model output.


\subsubsection{Elementary Effects}\label{sec:Morris}

The elementary effect test is a screening method able to overcome the main limitation of one-at-a-time designs and derivative-based methods, by considering wide ranges of variations for the inputs and averaging over a number of local derivative approximations, which allows to identify those input factors that are (i) negligible, (ii) linear and additive, or (iii) nonlinear or involved in interactions with other factors, with a relatively small number of model evaluations. 
Without loss of generality, we assume that the parameter domain is the $p$-dimensional unit hypercube. Hence, given a partition $\Omega_\iota$ of $[0,1]^p$ into $\iota\in\mathbb{N}$ levels, the elementary effect associated with the $i$-th input factor, for $i = 1,\dots,p$, is defined as
\begin{equation*}
    EE_i(\boldsymbol{\mu}) = \frac{y(\boldsymbol{\mu}+\mathbf{e}_i\Delta^\mu) - y(\boldsymbol{\mu})}{\Delta^\mu},
\end{equation*}
where $\boldsymbol{\mu}, \boldsymbol{\mu}+\mathbf{e}_i\Delta^\mu\in\Omega_\iota$, being $\mathbf{e}_i$ the unit vector of all zeros but the $i$-th component, and $\Delta$ the step size with value in $\left\{ 1/(\iota-1), \dots, (\iota-2)/(\iota-1)\right\}$. Since each elementary effect quantifies a local behavior, statistics of their distribution are estimated to obtain a global sensitivity measure. The mean and standard deviation of the distribution of $EE_i$ is constructed by randomly sampling different $\boldsymbol{\mu}$ from $\Omega_\iota$. Moreover, to avoid cancellation effects, the mean of the distribution $|EE_i|$ of the absolute values is also computed. Thus, for $r\in\mathbb{N}$ trajectories, that is, sets of points $\boldsymbol{\mu}_1,\dots,\boldsymbol{\mu}_r$ sampled in $\Omega_\iota$, we obtain the following sensitivity measures for the $i$-th input
\[
    m_i = \frac{1}{r}\sum_{j=1}^r EE_i(\boldsymbol{\mu}_j), 
\quad 
    m_i^* = \frac{1}{r}\sum_{j=1}^r 
\left\lvert EE_i(\boldsymbol{\mu}_j)\right\rvert,    
\quad
    sd_i = \frac{1}{r-1}\sum_{j=1}^r \left(EE_i(\boldsymbol{\mu}_j) - m_i\right)^2, \ \quad  i=1,\dots,p.
\]
In our setting, we fix the number of levels $\iota =6$ and choose $\Delta^\mu = \frac{\iota}{2(\iota-1)}$   to guarantee an equal probability sampling. Moreover, we consider $r=20$ trajectories, corresponding to $r(p+1)$ input-output evaluations, yielding a good balance between reliability of the results and computational resources required for the analysis.

Both $m_i$ and $m_i^*$ can be used to quantify the individual influence of the $i$-th input factor on the output $y$, whereas $sd_i$ estimates the ensemble of its effects due to nonlinearities and interactions with the other factors. Low values for $sd_i$ can be interpreted as the fact that the effect of the $i$-th parameter is independent of the point in which $EE_i$ is evaluated, and thus of the values taken by the other inputs.

\begin{figure}[h!]
    \centering
        \includegraphics[width=0.95\textwidth]{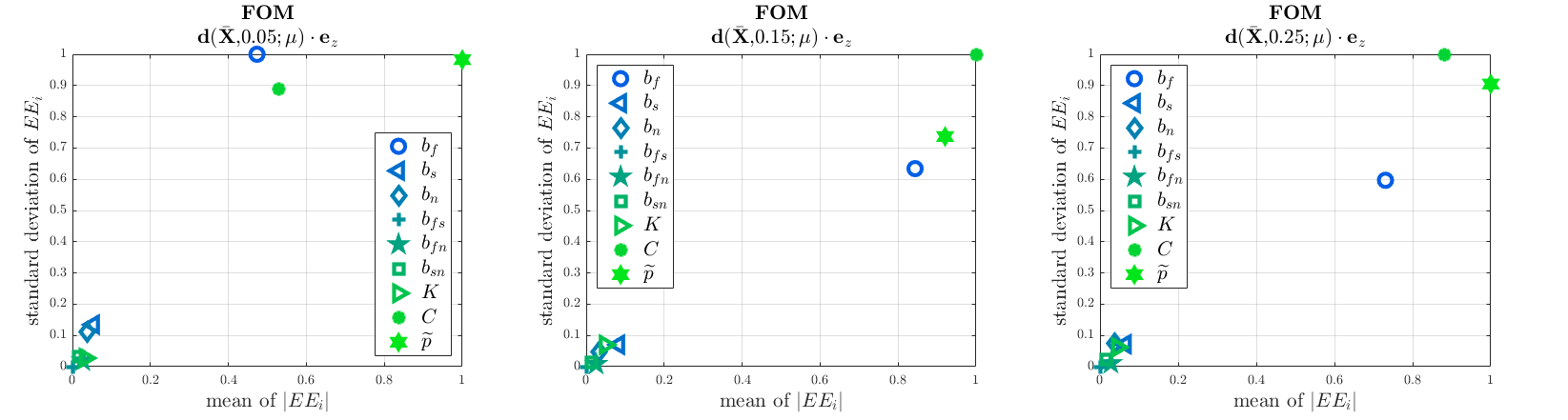}\\
    \includegraphics[width=\textwidth]{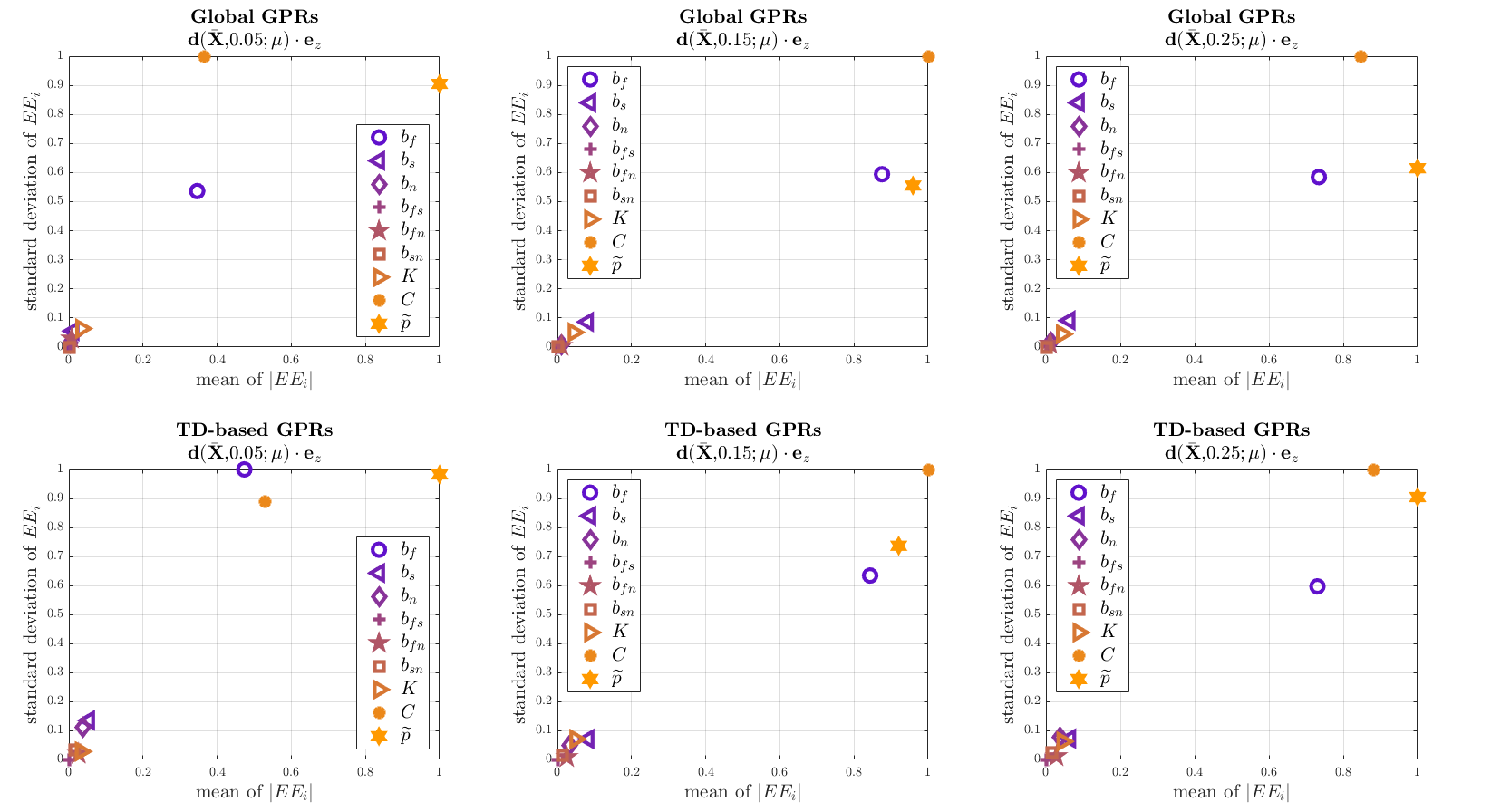}
    \caption{Deformation of a beam. Morris metrics computed using the outputs of the FOM (top), the global POD-GPR ROM (middle), and the tensor-decomposition-based  POD-GPR ROMs (bottom). The QoIs are given by the displacement along the $z$-axis of a point $P$, at three different time instants (from left to right, $t=0.05, 0.15, 0.25$).}
    \label{fig:beam_Morris}
\end{figure}

\paragraph{Test 1: deformation of a beam.} For the first problem, we report in Figure~\ref{fig:beam_Morris} (top panel) the estimated mean $\mathbf{m}^* = [m_1^*,\dots,m_p^*]^T$ and standard deviation $\mathbf{sd}^* = [sd_1,\dots,sd_p]^T$  of the elementary effects of each QoI computed using the FOM, whereas on the middle and bottom panels of Figure~\ref{fig:beam_Morris} we show the same statistical quantities computed using the POD-GPR ROMs outputs. To better visualize the results, a min-max scaling is applied to both $\mathbf{m}^*$ and $\mathbf{sd}$, for each model. Further details are reported in Table~\ref{tab:beam_morris}.

The elementary effects computed by means of tensor-decomposition-based approach are in great accordance with the high-fidelity ones. In any case, the ranking of the input parameters is always correctly identified. Relying on surrogates models allows to obtain accurate results of the SA in seconds, rather than hours, and thus reduce the size of the parameter space by identifying the less influential inputs.

\begin{table}[h!]
    \centering
    \begin{tabular}{|l|ccc|}
        \hline
        & FOM & Global GPR & TD-based GPR \\
        \hline
        Online CPU time & $3$ h $11$ min & $18$ s & $25$ s \\
        Speed-up & -- & $645$ & $450$ \\
        \hline
        mean of $|EE_1|$ ($b_f$) at $t=0.25$~s  & $2.5\cdot10^{3}$ & $2.0\cdot10^{3}$ & $2.3\cdot10^{3}$ \\
        mean of $|EE_8|$ ($C$) at $t=0.25$~s & $2.8\cdot10^{3}$ & $2.3\cdot10^{3}$ & $2.8\cdot10^{3}$ \\
        mean of $|EE_9|$ ($\widetilde{p}$) at $t=0.25$~s & $3.0\cdot10^{3}$ & $2.7\cdot10^{3}$ & $3.2\cdot10^{3}$ \\
        \hline
    \end{tabular}
        \caption{Deformation of a beam. Computational performance of the POD-GPR-ROM methods for the computation of the  elementary effects.}
    \label{tab:beam_morris}
\end{table}

\paragraph{Test 2: active contraction of a truncated ellipsoid.}	Due to the higher complexity of this  benchmark, we compute the elementary effects relying only on the POD-GPR ROMs. In particular, the global approach requires $18$~s to run the analysis, whereas the tensor-decomposition-based approach takes around $6$~min, almost $100\%$ of the time required to compute $r(p+1)=260$ input-output evaluations.

\begin{figure}[h!]
    \centering
    \includegraphics[width=\textwidth]{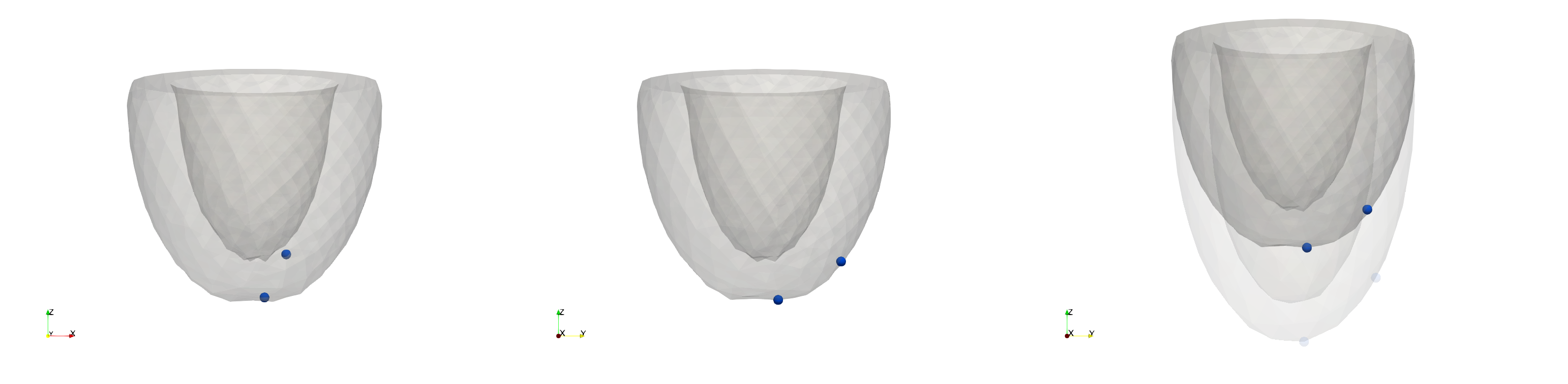}
    \caption{Active contraction of a truncated ellipsoid. Points $P_1$ and $P_2$ at time 0.25s from three different points of view.}
    \label{fig:prolate_QoI_selected}
\end{figure}

\begin{figure}[h!]
    \centering
    \includegraphics[width=\textwidth]{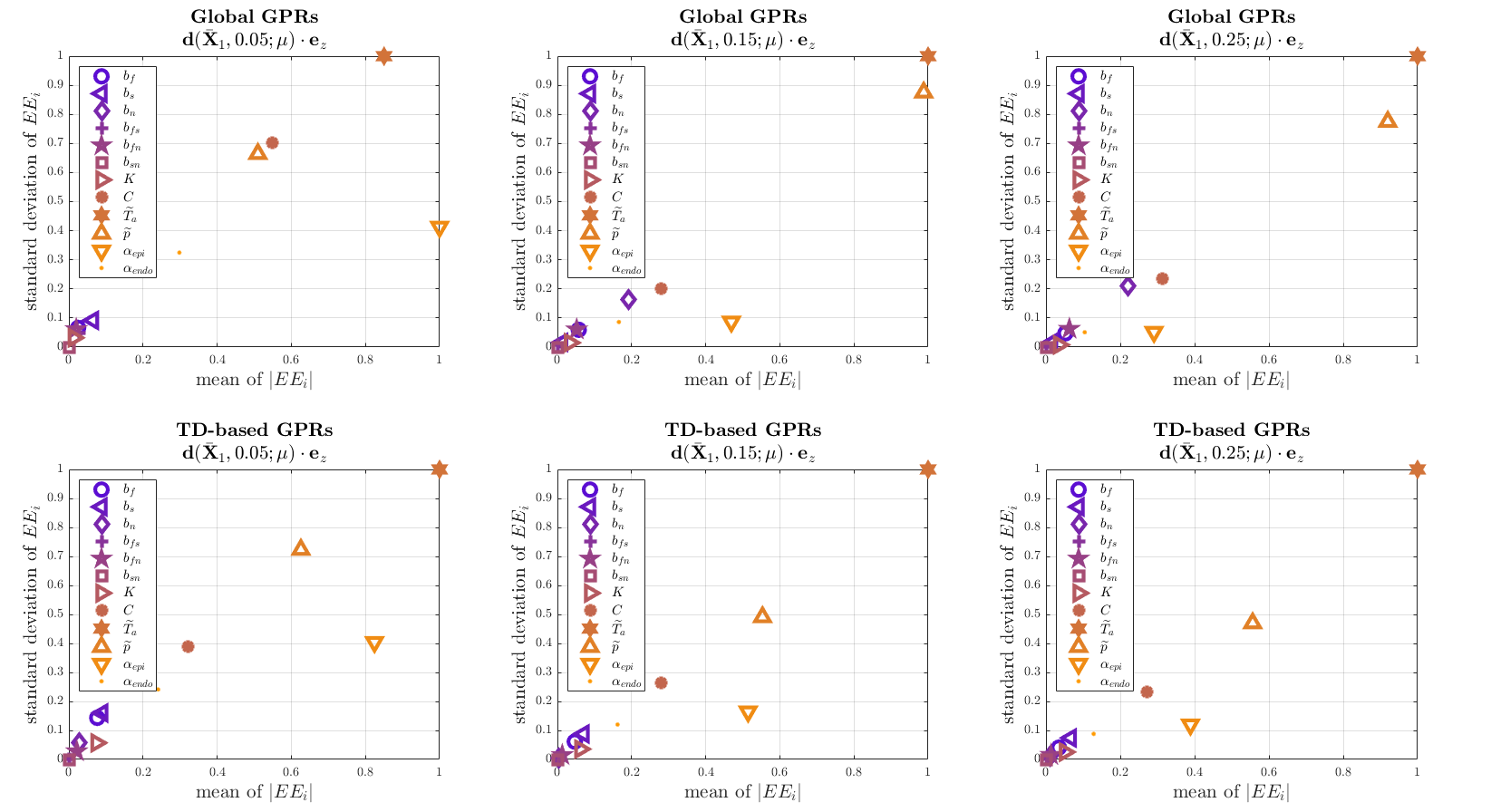}
    \caption{Active contraction of a truncated ellipsoid. Scatter plots of Morris metrics computed using the global (top) and the tensor-decomposition-based (bottom) POD-GPR ROMs. The QoIs are given by the displacement along the $z$-axis of a point $P_1$ at the epicardium located near the apex, at three different time instants.}
    \label{fig:prolate_Morris_apex}
\end{figure}

\begin{figure}[h!]
    \centering
    \includegraphics[width=\textwidth]{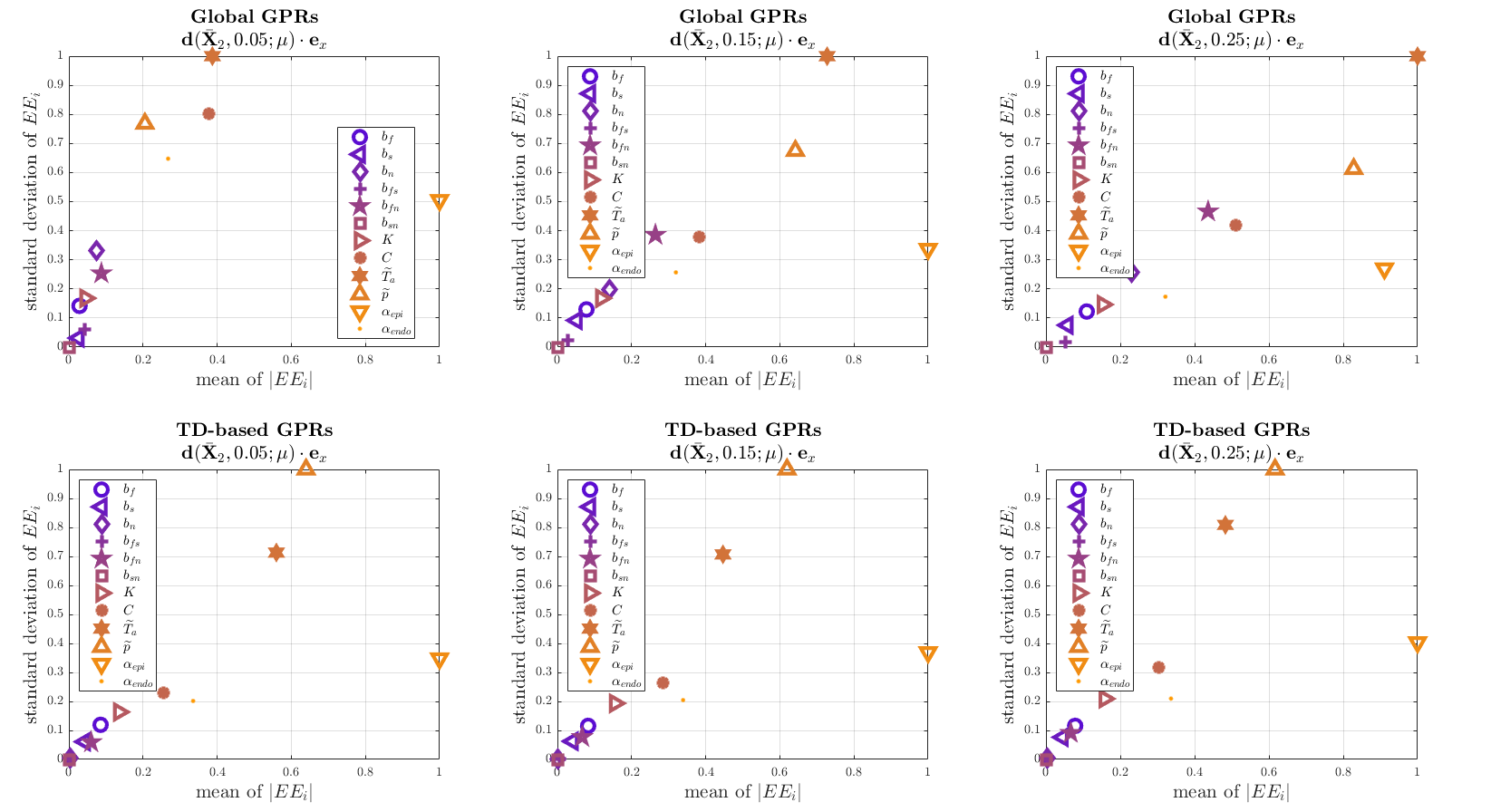}
    \caption{Active contraction of a truncated ellipsoid. Scatter plots of Morris metrics computed using the global (top) and the tensor-decomposition-based (bottom) POD-GPR ROMs. The QoIs are given by the displacement along the $x$-axis of a point $P_2$ at the epicardium located between the apex and the base, at three different time instants.}
    \label{fig:prolate_Morris_point}
\end{figure}

In Figures~\ref{fig:prolate_Morris_apex} and \ref{fig:prolate_Morris_point} we report the scatter plots of Morris indices computed using the POD-GPR ROMs for the QoIs associated with two different points located on the epicardium and highlighted in Figure~\ref{fig:prolate_QoI_selected} (namely, $P_1$ near the apex, and $P_2$ between the apex and the base, respectively). In particular, for the first point the vertical displacement at different time instances is considered, whether for the second point we focus on the displacement along the $x$-axis. The choice of the outputs of interest, differing in time, location or direction of the displacement, may have a great impact on the ranking of the input parameters, so that they have to be carefully chosen. However, by relying on POD-GPR approaches, one does not need to define the outputs in advance, since the training of the GPR is independent of the QoIs. Despite few discrepancies between the global and the tensor-decomposition-based approaches, we observe that the slope $\widetilde{T}_a$ of the active tension has large values for both the statistics, especially when taking into account the displacement along the $z$-axis. On the other hand, the fiber angle $\boldsymbol{\alpha}_{epi}$ shows a stronger influence on QoIs associated the the horizontal displacement. In any case, the less influential parameters are correctly located by both POD-GPR ROMs.

\subsubsection{Sobol' indices}

Variance-based methods are even more powerful tool to quantify the relative importance of individual factors, and are based on a decomposition of the variance of the model output $y$ into terms related to each input and to the interactions between them, that is
\begin{equation}\label{eq:ANOVA-HDMR}
    \text{Var}(y) = \sum_i V_i + \sum_{i<j} V_{ij} + \dots + V_{1\dots p},
\end{equation}
where $V_i = \text{Var}(\mathbb{E}[y|\boldsymbol{\mu}^{(i)}])$,  $V_{ij} = \text{Var}(\mathbb{E}[y|\boldsymbol{\mu}^{(i)},\boldsymbol{\mu}^{(j)}])-V_i-V_j$ , 
and so on, being $\boldsymbol{\mu}^{(i)}$ the $i$-th component of the parameter vector, for $i=1,\dots,p$. A primary measure of sensitivity is given by the \textit{first Sobol' index}
\begin{equation*}
    S_i = \frac{V_i}{\text{Var}(y)} = \frac{\text{Var}(\mathbb{E}[y|\boldsymbol{\mu}^{(i)}])}{\text{Var}(y)}, \quad i=1,\dots,p,
\end{equation*}
which measures the effect of varying the $i$-th input alone averaged over variations in all input parameters: the higher is $S_i$, the greater is the influence of $\boldsymbol{\mu}^{(i)}$ on the output $y$. Note that $\sum_i S_{i} \leq 1$. However, it is often impossible to separate the effects of the inputs on the output and one should look for higher order interactions. Rather than performing a full sensitivity analysis, which requires the evaluation of $2^p-1$ indices, we rely on the \textit{total effect indices}
\begin{equation*}
    S_{T_i} = \frac{\mathbb{E}[\text{Var}(y|\boldsymbol{\mu}^{(\sim i)})]}{\text{Var}(y)}
    = \frac{\text{Var}(y)-\text{Var}(\mathbb{E}[y|\boldsymbol{\mu}^{(\sim i)}])}{\text{Var}(y)}
    = 1 - \frac{\text{Var}(\mathbb{E}[y|\boldsymbol{\mu}^{(\sim i)}])}{\text{Var}(y)},
\end{equation*}
for $i=1,\dots,p$, where  $\boldsymbol{\mu}^{(\sim i)} = (\boldsymbol{\mu}^{(1)},\dots,\boldsymbol{\mu}^{(i-1)},\boldsymbol{\mu}^{(i+1)},\dots,\boldsymbol{\mu}^{(k)})$ denotes the random vector of all input factors but $\boldsymbol{\mu}^{(i)}$. In this case, the smaller is $S_{T_i}$, the less influential is $\boldsymbol{\mu}^{(i)}$, which can thus be arbitrarily fixed within its range of uncertainty without appreciably affecting the output of interest. Unlike $S_i$, it holds $\sum_i S_{T_i} \geq 1$, since $S_{T_{i_1}}$ and  $S_{T_{i_2}}$ ($i_1<i_2$) take both into account the interactions between $\boldsymbol{\mu}^{(i_1)}$ and $\boldsymbol{\mu}^{(i_2)}$.

In this work we employ the \textit{Saltelli method} \cite{saltelli2002making}, which is one of the most used and efficient procedure to evaluate Sobol' indices relying on $N_{samples}$ evaluations of the output quantity of interest, yielding a total number of $N_{samples}(p+2)$ output evaluations to compute both first Sobol' and total effect indices -- therefore entailing a huge gain compared to a crude Monte Carlo estimate of the variance-based sensitivity indices.

\begin{figure}[h!]
    \centering
    \includegraphics[width=0.95\textwidth]{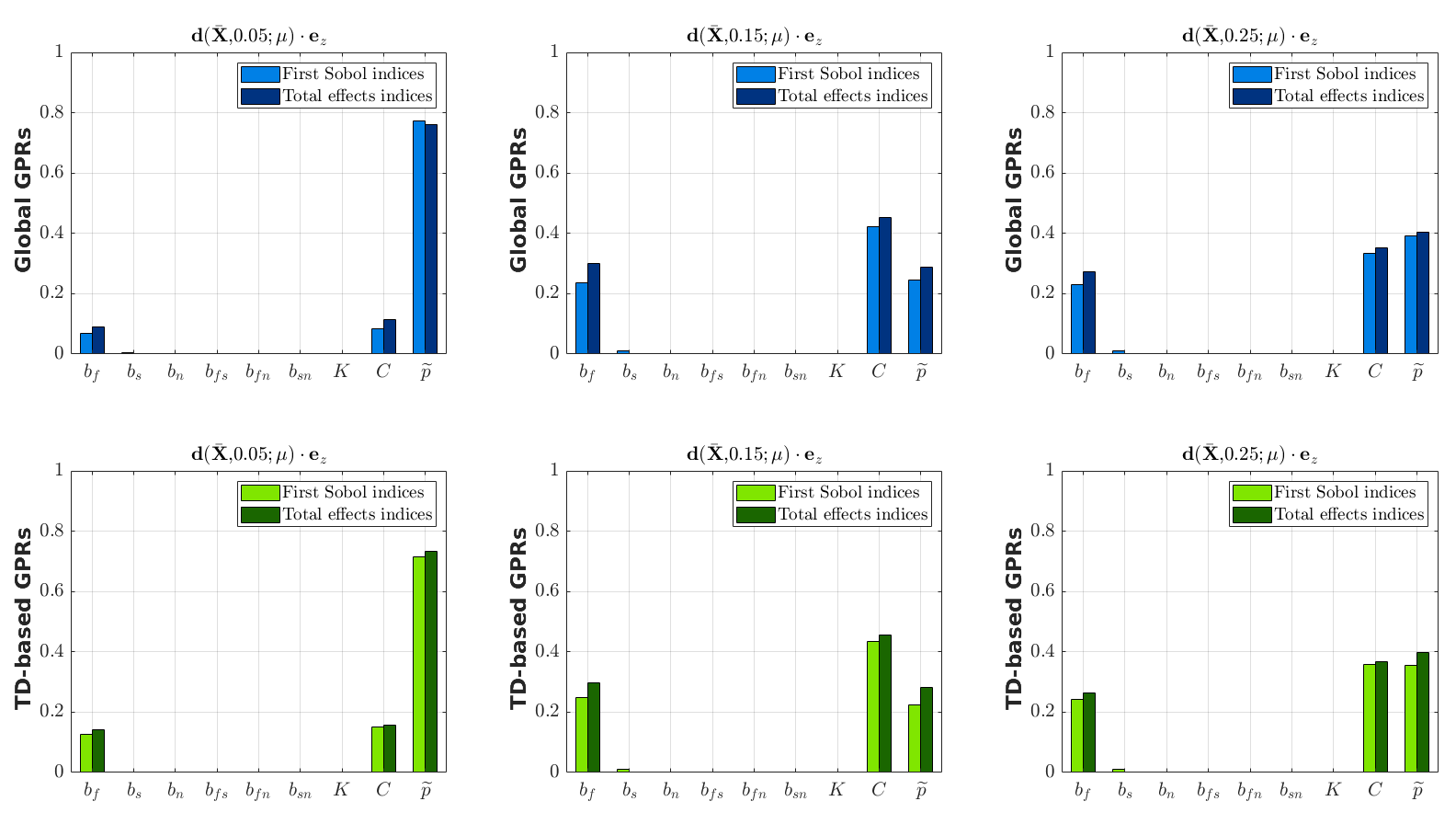}
    \caption{Deformation of a beam. Bar plot of Sobol' indices computed using the global (top) and the tensor-decomposition-based (bottom) POD-GPR ROMs. The QoIs are given by the displacement along the $z$-axis of a point $P$, at three different time instants.}
    \label{fig:beam_Sobol'}
\end{figure}

\paragraph{Test 1: deformation of a beam} We consider $N_{samples} = 200$ output evaluations, thus corresponding to $N_{samples}(p+2)=2200$ input-output evaluations, and rely on the global and the tensor-decomposition-based GPR models. Figure~\ref{fig:beam_Sobol'} confirms the results of the SA previously conducted using Morris elementary effects, that is the scaling factor $C$, the slope $\widetilde{p}$ of the pressure load and the material stiffness in the fiber direction $b_f$ are the most influential parameters for the test at hand, whereas little to no influence on the selected output is given by the other parameters, which can thus be fixed in their domain. Further details are reported in Table~\ref{tab:beam_Sobol}.

 \begin{table}[h!]
    \centering
    \begin{tabular}{|l|ccc|}
        \hline
        & FOM & Global GPR & TD-based GPR \\
        \hline
        Online CPU time & $1$ d $8$ h $45$ min & $147$ s & $290$ s \\
        Speed-up & -- & $804$ & $406$ \\
        \hline
        $S_{1}$ ($b_f$) at $t=0.25$~s  & $0.25$ & $0.23$ & $0.24$ \\
        $S_{8}$ ($C$) at $t=0.25$~s & $0.37$ & $0.33$ & $0.36$ \\
        $S_{9}$ ($\widetilde{p}$) at $t=0.25$~s & $0.35$ & $0.39$ & $0.36$ \\
        \hline
        $S_{T_1}$ ($b_f$) at $t=0.25$~s  & $0.28$ & $0.27$ & $0.26$ \\
        $S_{T_8}$ ($C$) at $t=0.25$~s & $0.38$ & $0.35$ & $0.37$ \\
        $S_{T_9}$ ($\widetilde{p}$) at $t=0.25$~s & $0.38$ & $0.41$ & $0.40$ \\
        \hline
    \end{tabular}
        \caption{Deformation of a beam. Computational performance of the POD-GPR-ROM methods for the computation of Sobol' indices.}
    \label{tab:beam_Sobol}
\end{table}

\begin{figure}[h!]
    \centering
    \includegraphics[width=0.95\textwidth]{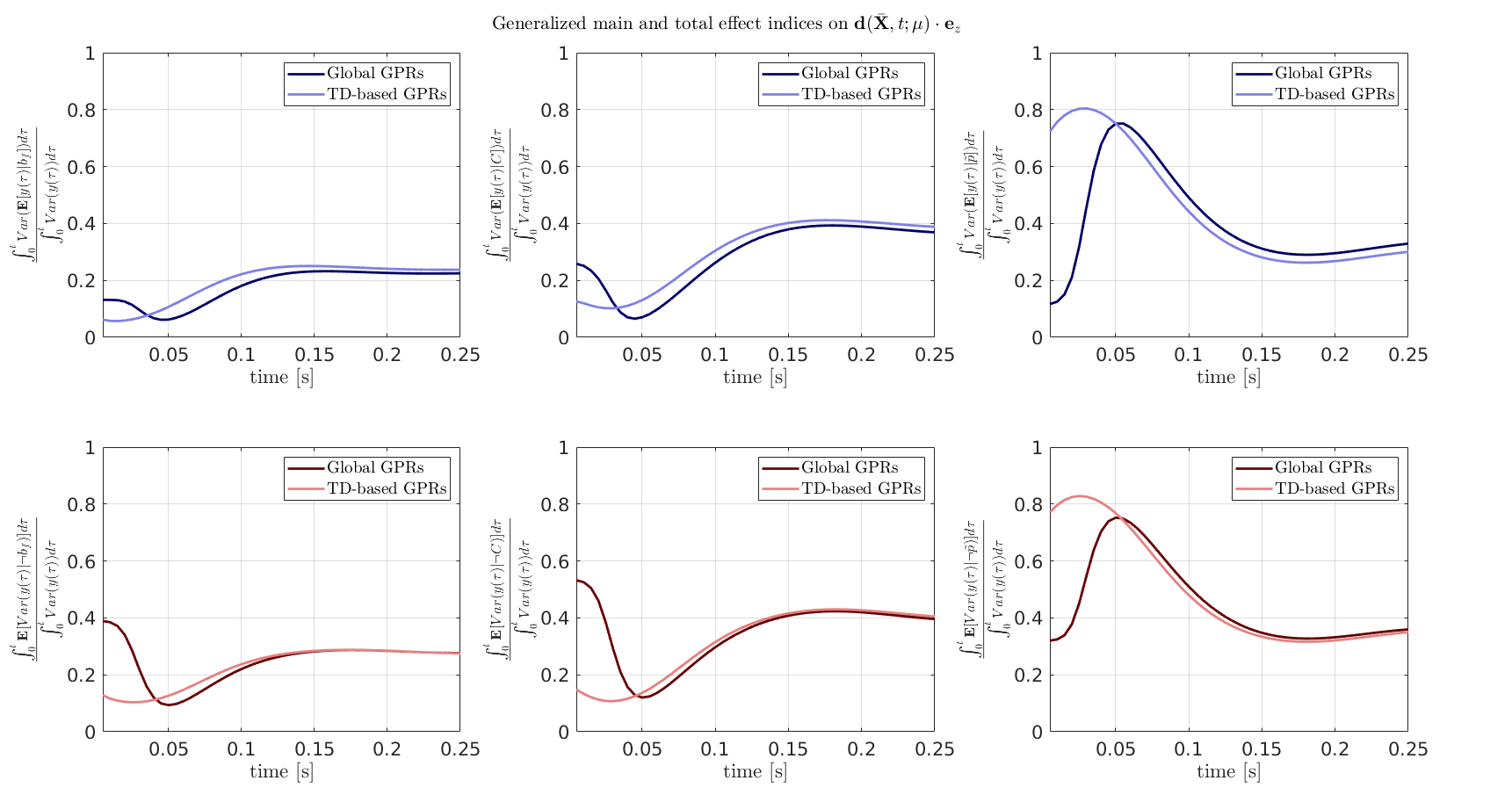}
    \caption{Deformation of a beam. Evolution in time of Sobol' main (top) and total (bottom) effect indices for $b_f$, $C$ and $\widetilde{p}$ computed using the global and the tensor-decomposition-based POD-GPR ROMs.}
    \label{fig:beam_Sobol_generalized}
\end{figure}

Furthermore, we compute Sobol' indices for all $n=1,\dots,N_t$. Figure~\ref{fig:beam_Sobol_generalized} shows how the following generalization of Sobol' indices
\begin{equation*}
    \frac{\int_0^t \text{Var}(\mathbb{E}[y(\tau)|\boldsymbol{\mu}^{(i)}]) d\tau}{\int_0^t \text{Var}(y(\tau)) d\tau}, \quad \frac{\int_0^t\mathbb{E}[\text{Var}(y(\tau)|\boldsymbol{\mu}^{(\sim i)})]d\tau}{\int_0^t\text{Var}(y(\tau))},
\end{equation*}
evolve over time for the most effective inputs, i.e., $b_f$, $C$ and $\widetilde{p}$. We observe that the indices associated with the material inputs increase as time grows, whereas the pressure parameter reduces its influence on the output as the maximum slope is reached. However, when interested in, e.g., time-dependent model outputs, performing independent SAs on each prediction associated with a given time step $t^n$, with $n=1,\dots,N_t$, may lead to unsatisfactory results. In fact, the $N_t$ variables are correlated, so that important information cannot be extracted from separated analysis. In this scenario, one needs to rely on methods specific for SA on functional or multivariate
outputs, such as \cite{campbell20061468}, where the output is expanded in a proper basis and  standard SA methods are then applied to the coefficients of this expansion, and \cite{gamboa2014sensitivity, lamboni2011multivariate, alexanderian2020variance}, where generalized Sobol' sensitivity indices are proposed. 

\begin{figure}[h!]
    \centering
    \includegraphics[width=0.95\textwidth]{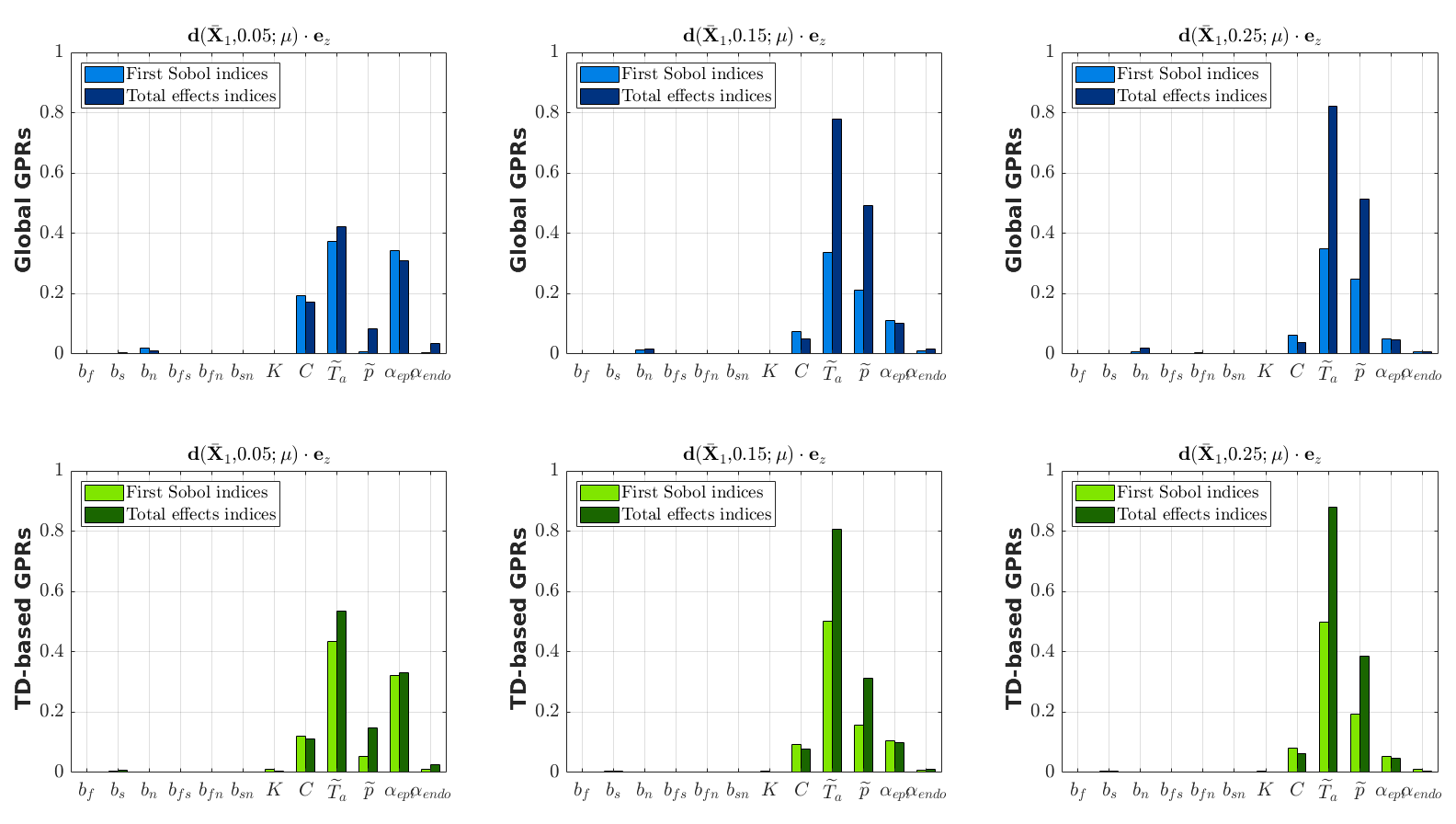}
    \caption{Active contraction of a truncated ellipsoid. Bar plot of Sobol' indices computed using the global (top) and the tensor-decomposition-based (bottom) POD-GPR ROMs. The QoIs are given by the displacement along the $z$-axis of a point $P_1$ at the epicardium located near the apex, at three different time instants.}
    \label{fig:prolate_Sobol_apex}
\end{figure}

\begin{figure}[h!]
    \centering
    \includegraphics[width=0.95\textwidth]{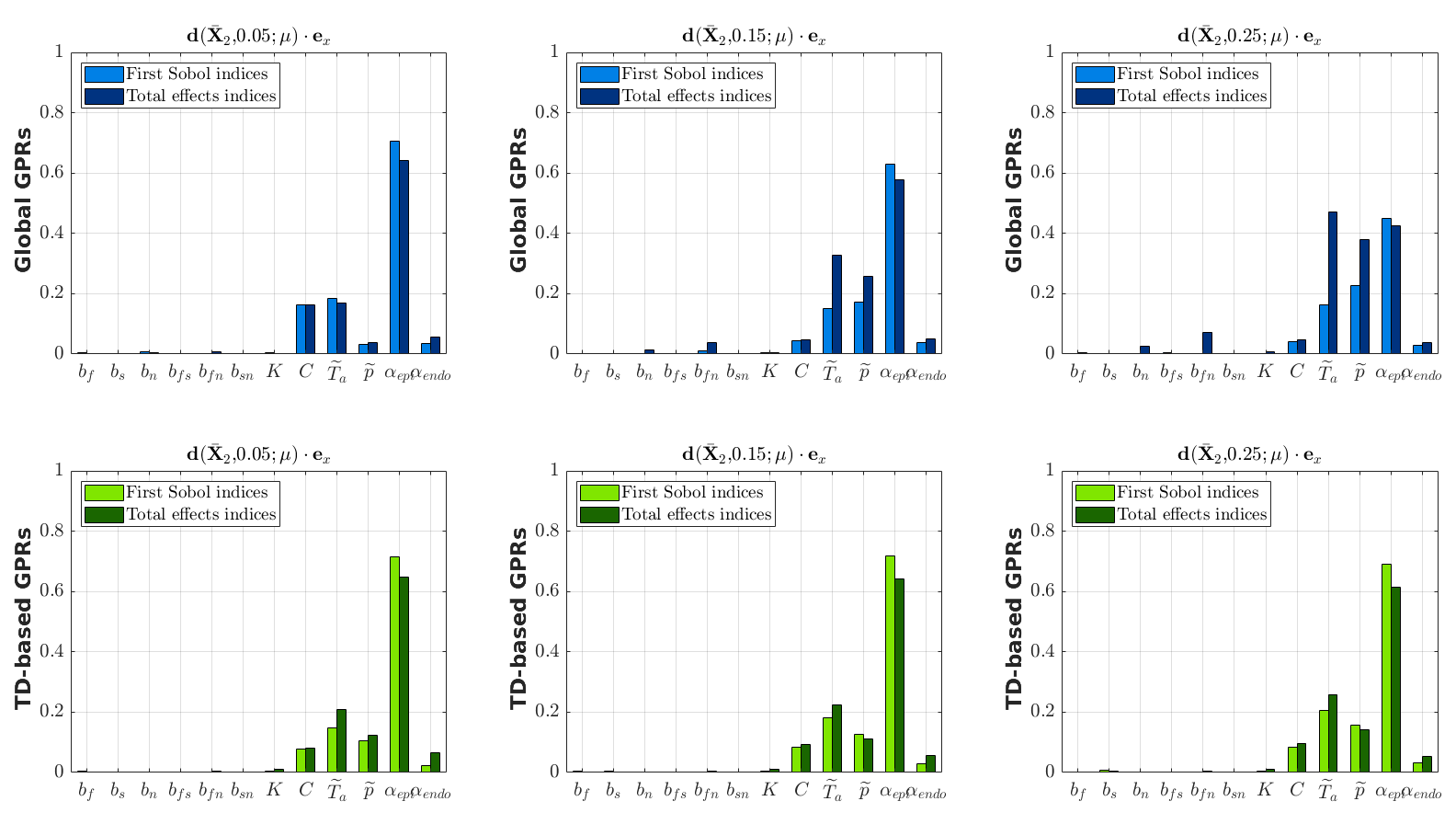}
    \caption{Active contraction of a truncated ellipsoid. Bar plot of Sobol' indices computed using the global (top) and the tensor-decomposition-based (bottom) POD-GPR ROMs. The QoIs are given by the displacement along the $x$-axis of a point $P_2$ at the epicardium located between the apex and the base, at three different time instants.}
    \label{fig:prolate_Sobol_point}
\end{figure}

\paragraph{Test 2: active contraction of a truncated ellipsoid} For the second benchmark we consider again $N_{samples}=200$, corresponding in this case to $N_{samples}(p+2)=2800$ runs of the POD-GPR ROMs. The results obtained with Sobol SA regarding the ranking of the input parameters, especially the detection of the less influential ones, is in agreement with Morris screening. Figures~\ref{fig:prolate_Sobol_apex} and \ref{fig:prolate_Sobol_point} show the sensitivity indices obtained for the QoIs related to the same points reported in Section~\ref{sec:Morris}, showing that $\widetilde{T}_a$ has a great influence, in particular on the displacement towards the base; moreover, it exhibits strong nonlinear interactions, as pointed out by the large value of the total effects at different times. Regarding the fiber angle an the epicardium $\boldsymbol{\alpha}_{epi}$, its impact is clearly visible on the horizontal displacement of the second point, which is located between the apex and the base.  Almost all the parameters related to the material law have very small indices, meaning that they can be arbitrarily fixed withing their range of variation without significantly affect the computation of the outputs of interest.

\subsection{Parameter estimation}\label{sec:inverseUQ}

Finally, we address the solution of parameter estimation problems in a Bayesian framework, thus describing all inputs as random variables and seeking their probability density functions (PDFs), which, starting from a prior intuition, are updated as new data becomes available. 

Given an observation $\mathbf{y}_{obs}$, corresponding to the realization of the output random variable $\mathbf{y}$, we want to infer the underlying probability distribution that produced the data. Note that, with respect to the previous section, we do not assume $\mathbf{y}$ to be a univariate model output, hence the change in notation. In particular, we consider $\mathbf{y}_{obs}$ to be the vector of QoIs evaluated on the FOM solution $\mathbf{u}_h(\boldsymbol{\mu}^*)$, computed for a specified parameter vector $\boldsymbol{\mu}^*\in\mathcal{P}$. To mimic experimental error, we consider an additive independent identically distributed Gaussian noise $\boldsymbol{\varepsilon}\sim\mathcal{N}(0,\boldsymbol{\Sigma}_{\boldsymbol{\varepsilon}})$, so that $\mathbf{y}_{obs} = \mathbf{y}(\boldsymbol{\mu}^*) + \boldsymbol{\varepsilon}$. The \textit{posterior} PDF of the random vector $\boldsymbol{\mu}$ can be obtained using Bayes' formula as
\begin{equation*}
    \pi(\boldsymbol{\mu}|\mathbf{y}_{obs}) \propto \pi(\mathbf{y}_{obs}|\boldsymbol{\mu})\pi_0(\boldsymbol{\mu}),
\end{equation*}
where
\begin{itemize}
    \item $\pi(\mathbf{y}_{obs}|\boldsymbol{\mu})$ is the \textit{likelihood};
    \item $\pi_0(\boldsymbol{\mu})$ is the \textit{prior} PDF of the inputs, which reflects any previous knowledge one might have about the parameters. If no previous knowledge is available, non-informative priors can be used.
\end{itemize}
Thanks to the assumption on the measurement error, the likelihood is given by $\pi(\mathbf{y}|\boldsymbol{\mu}) = \pi_{\boldsymbol{\varepsilon}}(\mathbf{y}-y(\boldsymbol{\mu}))$, where $\pi_{\boldsymbol{\varepsilon}}$ is the PDF of the noise $\boldsymbol{\varepsilon}$. Finally, one can compute useful statistical indicators, such as the posterior mean and covariance, by repeatedly sampling from the posterior distribution.

Markov Chains Monte Carlo (MCMC) methods represent a well-known class of techniques used for sampling from a probability distribution. In this work we rely on the Metropolis-Hastings algorithm, whose key idea is to generate a sequence of $N_{MC}$ samples (with $N_{MC}$ fixed, but sufficiently large) such that, at each iteration, the new candidate value is chosen according to a specified proposal distribution $\pi_{prop}$, and it is either accepted or rejected with some probability determined by the properties of the likelihood and the prior. To reduce the bias introduced by the choice of the initial sample, we perform \textit{burn-in}, i.e., we discard $N_{burn\text{-}in}$ iterations at the beginning of an MCMC run; moreover, we keep only one sample every $N_{thin}$, so that the total number of chain iterations is given by $\lfloor(N_{MC}-N_{burn\text{-}in})/N_{thin}\rfloor$.

\paragraph{Test 1: deformation of a beam.} Going back to the problem of the deformation of a beam discussed in Section~\ref{sec:beam}, we employ the global and the tensor-decomposition-based POD-GPR ROMs built for $N_s=50$, which are characterized by a time-average $L^2(\Omega_0)$-absolute error of $4\cdot10^{-4}$ with respect to the FOM for the solution to the forward problem. The values of the target parameter $\boldsymbol{\mu}^*$ used to generate the observation $\mathbf{y}_{obs}$ and those of the chain starting point $\boldsymbol{\mu}^0$ are reported in Table~\ref{tab:beam_MCMC_points}. The inputs to the Metropolis-Hastings algorithm are instead listed in Table~\ref{tab:beam_MH}.

\begin{table}[h!]
    \centering
    \begin{tabular}{|c|ccccccccc|}
        \hline
        & $b_f$ & $b_s$ & $b_n$ & $b_{fs}$ & $b_{fn}$ & $b_{sn}$ & $K$ & $C$ & $\widetilde{p}$ \\
        & & & & & & & [kPa] & [kPa] & [kPa]  \\
        \hline
        $(\cdot)^*$ & $6.2$ & $1.2$ & $2.8$ & $5.8$ & $5.8$ & $2.8$ & $27$ & $1.2$ & $0.0058$ \\
        $(\cdot)^0$ & $8$ & $2$ & $2$ & $4$ & $4$ & $2$ & $50$ & $2$ & $0.004$ \\
        \hline
    \end{tabular}
    \caption{Deformation of a beam. Target and starting values of the Markov chain for the input parameters.}
    \label{tab:beam_MCMC_points}
\end{table}

\begin{table}[h!]
    \centering
    \begin{tabular}{|l|r||l|l|} 
        \hline
        $N_{MC}$ & 10000& $\pi_0$ & $\mathcal{U}$[$\mathcal{P}$] (non-informative prior)\\
        \hline
        $N_{burn\text{-}in}$ & 500 & $\pi_{\boldsymbol{\varepsilon}}$ & $\mathcal{N}(0,\sigma_{\boldsymbol{\varepsilon}}^2\mathbb{I})$, with $\sigma_{\boldsymbol{\varepsilon}}^2 = 10^{-5}$ or $10^{-6}$\\
        \hline
        $N_{thin}$ & 4 & $\pi_{prop}$ & $\mathcal{U}$[$\mathcal{P}$]\\
        \hline		
    \end{tabular}
    \caption{Deformation of a beam. Inputs to the Metropolis-Hastings algorithm.}
    \label{tab:beam_MH} 
\end{table}
	
\begin{figure}[h!]
    \centering
    \includegraphics[width=\textwidth]{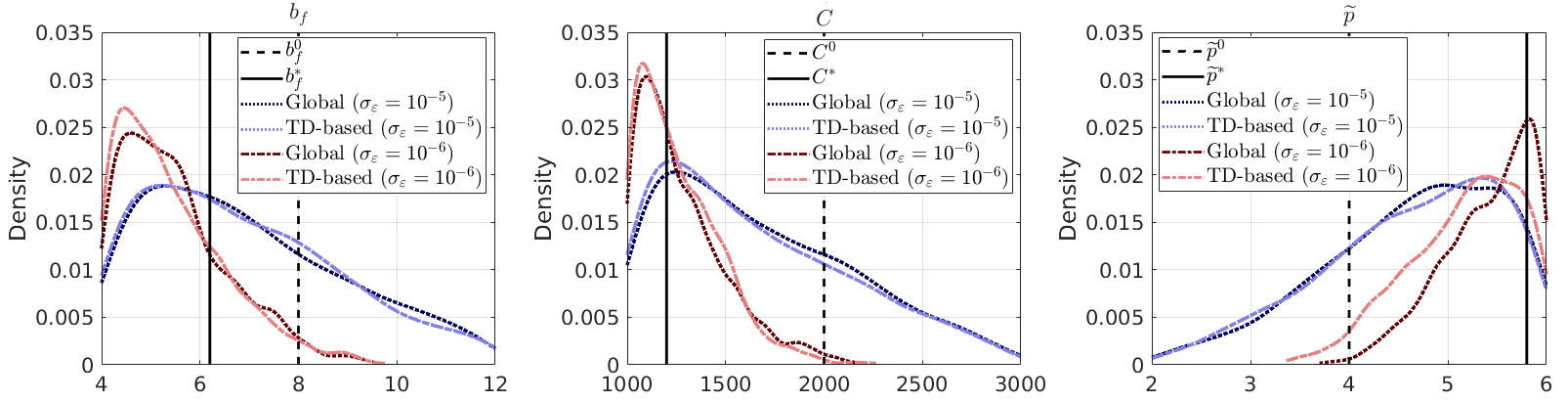}
    \caption{Deformation of a beam. Fitted posterior density function obtained from the MCMC samples of the most influential input parameters computed using the global and tensor-decomposition-based POD-GPR ROMs, for two different values of $\sigma_\varepsilon$. The vertical lines shows the chain starting point (dashed) and the target values (solid).}
    \label{fig:beam_MCMC_pdf}
\end{figure}

The samples of the MCMC chain obtained after the burn-in and the thinning operations are used to compute the posterior density functions reported in 
Figure~\ref{fig:beam_MCMC_pdf} for $\boldsymbol{\Sigma}_{\boldsymbol{\varepsilon}}=10^{-5}\mathbb{I}$ and $10^{-6}\mathbb{I}$, where $\mathbb{I}\in\mathbb{R}^{N_{QoI}\times N_{QoI}}$ denotes the identity matrix and $N_{QoI}=3$. Moreover, the chain starting point $\boldsymbol{\mu}^0$ and the target parameter $\boldsymbol{\mu}^*$ are plotted as a black dashed line and a black solid line, respectively, for each parameter. We observe that, using the POD-GPR ROMs, we are able to correctly identify the parameter related to the material stiffness in the fiber direction $b_f$, the multiplicative constant $C$ and the slope of the pressure $\widetilde{p}$, whereas the other parameters are not captured. These results confirm those of the SA reported, e.g., in Figure~\ref{fig:beam_Sobol'}, which identify $b_f,C$ and $\widetilde{p}$ as the most influential inputs, while the others as non-influential. For what concerns the CPU time required for the generation of the MCMC chain, the global approach requires less that $12$~min of computation, while the tensor-decomposition-based approach has a total run time around $23$~min. Relying on the FOM would have required $8$~d $17$~h, showing how crucial is the need to rely on efficient reduced models in multi-query scenarios.

\paragraph{Test 2: active contraction of a truncated ellipsoid.}
To conclude, we perform parameter estimation in the case of a contracting truncated ellipsoid. However, due to the high computational complexity of this benchmark test and the large number of numerical simulations required, we rely only on the global POD-GPR based model, which is characterized by fast CPU times during the online stage. We recall that, in this case, the time-average $L^2(\Omega_0)$-absolute error is around $3\cdot10^{-2}$ with respect to the FOM.

The chain starting point and the target are reported in Table~\ref{tab:prolate_reference_values}. Moreover, assuming some previous knowledge about the input values that have generated the observation, we reduce the support of the uniform distribution with respect to the whole parameter domain. The numerical setting used for the MCMC simulation is the same as in the previous test case reported in Table~\ref{tab:beam_MH}, apart from the number of chain samples, which is increased to $N_{MC}=25000$, and the standard deviation of the Gaussian noise for which we now defined with $\sigma_\varepsilon=10^{-7}$.

\begin{table}[h!]
    \centering
    \begin{tabular}{|c|cccccccccccc|}
        \hline
        & $b_f$ & $b_s$ & $b_n$ & $b_{fs}$ & $b_{fn}$ & $b_{sn}$ & $K$ & $C$ & $\widetilde{p}$ & $\widetilde{T}_a$ & $\boldsymbol{\alpha}_{epi}$ & $\boldsymbol{\alpha}_{endo}$ \\
        & & & & & & & [kPa] & [kPa] & [kPa] & [kPa] & [$^\circ$]& [$^\circ$] \\
        \hline 
        $(\cdot)^*$ & $7.16$ & $2.21$ & $1.79$ & $4.42$ & $3.58$ & $2.21$ & $56$ & $1.7$ & $66.3$ & $14.4$ & $-99.3$ & $80.7$\\
        \hline
        $(\cdot)^0$ & $8$ & $2$ & $2$ & $4$ & $4$ & $2$ & $50$ & $2$ & $60$ & $15$ & $-90$ & $90$\\
        \hline
    \end{tabular}
    \caption{Active contraction of a truncated ellipsoid. Target and starting values of the Markov chain for the input parameters.}
    \label{tab:prolate_reference_values}
\end{table}
 
\begin{figure}[h!]
    \centering
    \includegraphics[width=\textwidth]{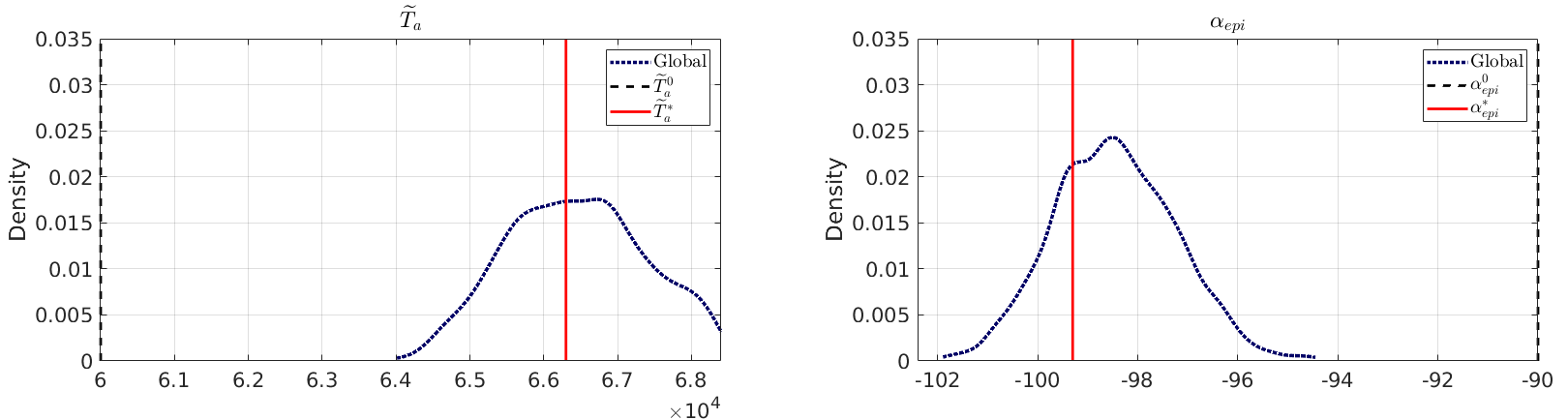}
    \caption{Active contraction of a truncated ellipsoid. Fitted posterior density function obtained from the MCMC samples of the most influential input parameters computed using the global POD-GPR ROM, for $\sigma_\varepsilon = 10^{-7}$. The vertical lines shows the chain starting point (dashed) and the target values (solid).}
    \label{fig:prolate_MCMC_1e-7_pdf}
\end{figure}

The PDF computed using the global GPR model, the chain starting point $\boldsymbol{\mu}^0$ and the target parameter $\boldsymbol{\mu}^*$ are reported in Figure~\ref{fig:prolate_MCMC_1e-7_pdf} using a dark 
blue dotted line, a black dashed line and a 
black solid line, respectively, for each input parameter. We observe that both the active tension $\widetilde{T}_a$ and the fiber angle at the epicardium $\boldsymbol{\alpha}_{epi}$ are correctly inferred, as we might expect from the SA results, which identified these inputs to be overall the most influential ones. Coherently, the identification of the remaining parameters is in general worse. For what concerns the CPU time required for the generation of the MCMC chain, the global approach requires $25$~min of computation, while relying on the FOM would have required more than $100$~days, becoming computationally unaffordable if the dimension $N_h$ of the underlying problem is further increased.

\section{Conclusions}\label{sec:conclusions}
In this work we have addressed the solution to parameterized, nonlinear, time-dependent problems by means of POD-GPR ROMs, able to capture with sufficient accuracy the state solution dynamics at much lower computational cost with respect to high-fidelity FOMs. In particular, exploiting GP regression to approximate the POD projection coefficient allows to avoid intrusive hyper-reduction techniques, although still providing a way to embed physics into the ROM. Another advantage of these methods is that, on one hand, they do not need to compute the whole solution dynamics when require only at specific points in space or time, and, on the other hand, one does not need to specify in advance the output quantities, being the training of the GPs independent of them, so that the QoIs can be chosen online according to the specific application.

Regarding the construction of the GPRs, we have compared a global approach and a tensor-decomposition-based one. Different kernels have been tested, as well as scaling techniques for the input/output data to the GPs, showing that a ARD-RBF covariance function and the standardization scaling are to be preferred. Finally, a convergence analysis with respect to the training set size has allowed to find a good trade-off between accuracy and online efficiency.

Both POD-GPR ROMs show low approximation errors with respect to the high-fidelity FOM, being able to accurately predict the system solution for new input values within the parameter space, as well as time instances not used during training. Regarding the computational costs, the tensor-decomposition-based approach is the most efficient during the offline stage, especially as the size of the problem increases. On the other hand, global POD-GPR ROMs shows higher online efficiency, and thus have to be preferred when hundreds or thousand input-output evaluations are required. Finally, the capabilities of the POD-GPR ROMs have been successfully tested on the solution to multi-query problems, namely global sensitivity analysis and parameter estimation.

\section*{Acknowledgment}
LC and PZ acknowledge the partial support from Regione Lombardia project NEWMED under Grant No. POR FESR 2014–2020. 
AM, SF, PZ are members of Gruppo Nazionale per il Calcolo Scientifico (GNCS) of Istituto Nazionale di Alta Matematica (INdAM). 
SF is supported by the PNRR-PE-AI FAIR project funded by the NextGeneration EU program.
AM and PZ acknowledge the partial support of the project `Sviluppo di sinergie fra Calcolo Scientifico e Machine Learning per applicazioni biomediche' funded by GNCS. 
MG acknowledges the financial support from \emph{Sectorplan Bèta} (the Netherlands) under the focus area \emph{Mathematics of Computational Science}.

\bibliographystyle{alpha}
\bibliography{bibliography}

\end{document}